\documentclass[a4,11pt]{article}

\font\elevensf=cmss10 scaled\magstephalf

\newtheorem{theorem}{{\elevensf THEOREM}}[section]

\newtheorem{lemma}{{\elevensf LEMMA}}[section]
\newtheorem{corollary}{{\elevensf COROLLARY}}[section]

\newtheorem{remark}{{\elevensf REMARK}}[section]

\newcommand{\qed}{\hfill \rule{1.6mm}{1.6mm}}

\def\CC{{\rm \kern.24em \vrule width.02em height1.4ex
depth-.05ex \kern-.26em C}}

\def\TagOnRight

\def\AA{{\it I}\hskip-3pt{\tt A}}

\def\QQ{\rlap {\raise 0.4ex \hbox{$\scriptscriptstyle |$}}
  {\hskip -0.1em Q}}

\newcommand{\be}{\begin{equation}}
\newcommand{\ee}{\end{equation}}
\newcommand{\bea}{\begin{eqnarray}}
\newcommand{\eea}{\end{eqnarray}}
\newcommand{\Bea}{\begin{eqnarray*}}
\newcommand{\Eea}{\end{eqnarray*}}

\newcommand{\de}{\delta}

\newcommand{\la}{\lambda}
\newcommand{\si}{\sigma}

\renewcommand{\th}{\theta}
\newcommand{\nn}{\nonumber}

\newcommand{\ep}{\epsilon}
\catcode`\@=11

\def\theequation{\@arabic{\c@section}.\@arabic{\c@equation}}

\catcode`\@=12
\title{
\vskip-3cm
  Asymptotic Properties of an Estimator of the Drift Coefficients of Multidimensional
  Ornstein-Uhlenbeck  Processes that are not Necessarily Stable
  }
  
\author{Gopal K. Basak\footnote{(Corresponding author) Stat-Math Unit,
Indian Statistical Institute, Kolkata 700~108, India. 
 E-mail:~gkb@isical.ac.in}
\ and \  
Philip Lee\footnote{ 
 JPMorganChase Bank, N.A.,
 Asia Rates Strategy, Hong Kong.
 Email: philip.pk.lee@jpmorgan.com}\ 
}
\date{}

\begin{document}

\maketitle

\parskip10pt
\parindent0in
\vspace*{-1.0cm}
\noindent 
\begin{abstract}
{
In this paper, we investigate the consistency 
and asymptotic efficiency
of an estimator
 of the drift matrix, $F$, of Ornstein-Uhlenbeck processes that are not
 necessarily
  stable. We consider all the cases. (1) The eigenvalues of $F$ are in 
the right half space (i.e., eigenvalues with positive real parts).
In this case the process grows exponentially fast.
    (2) The eigenvalues of $F$ are on the left half space 
(i.e., the eigenvalues with negative or zero real parts).
 The process where all eigenvalues of $F$ have negative real parts
is called a stable process and has a unique invariant (i.e., stationary) distribution.
 In this case the process does not grow.
 When the eigenvalues of $F$ have zero real parts (i.e., the case of zero
eigenvalues and purely imaginary eigenvalues) the process grows polynomially
fast.
Considering (1) and (2) separately,
we first show that an estimator, $\hat{F}$,
 of $F$ is consistent. We then combine them
to present 
 results for the
 general Ornstein-Uhlenbeck processes.
 We adopt similar procedure to show the asymptotic efficiency of the estimator.
}

\end{abstract}

\vspace*{-0.4cm}


\vspace*{-0.2cm}


\section{Introduction}
\setcounter{equation}{0}
\setcounter{section}{1}
\label{sec.1}
\vspace{-0.4cm}
Multidimensional 
processes with linear drift parameter have been used for modelling 
 various physical phenomena. 
 Among recent papers, works by
Jankunas and Khasminskii
(\cite{jnkh})
and  Khasminskii, Krylov and Moshchuk (\cite{khkm})
 on the estimation of the drift parameters of linear stochastic
differential equations 
(of the form,  $d X_t = A X_t dt + \sum_{i=1}^n 
\si_i X_t d w_i (t)$
and $d X_t = A_{\th} X_t dt + \sum_{i=1}^m 
\si_i X_t d w_i (t)$) can be mentioned.
It
 should be noted that
 our work on Ornstein-Uhlenbeck (OU) processes
 does not follow from theirs and that the methodology used in our paper
 is also quite different from theirs.

The motivation for this work comes from Lai and Wei's paper \cite{key1},
in which the authors
 have shown the strong consistency of the least square estimators of
the coefficients of the discrete univariate general AR(p) processes.
In this paper, we not only show that an estimator (which is 
the maximum likelihood estimator in the special case when $A$ is nonsingular)
 of the drift
 parameter
of the general multidimensional OU process is consistent but also show that it is asymptotically efficient.
We consider the following SDE representation of the OU process:
\begin{eqnarray}
\label{eq1.1}
dY_t = FY_t dt + A dW_t 
\end{eqnarray}
with any starting point $Y_0$ independent of the Brownian motion
 $\{W_t, \ \ t \geq 0 \}$. Here $Y$ is a $p$-dimensional process,
$A$ is a constant matrix of $p\times r$ dimesnion and
$W_t$ is a $r$-dimensional standard Brownian motion. 
Notice that it is always
easier to estimate $A$ through quadratic variation of the process by
 using It\^{o}'s rule. But, estimating $F$ is usually the
 more difficult
 task.
It is generally believed that one needs stationarity of the process
 to estimate $F$.
However, one may observe,
$\int_{0}^{T} dY_t Y_t' = F (\int_{0}^{T} Y_{t} Y_{t}' dt ) + 
A (\int_{0}^{T} dW_t Y_t')$.
 Thus, we define, 
$\hat{F}_T = (\int_{0}^{T} dY_t Y_t') (\int_{0}^{T} Y_{t} Y_{t}' dt )^{-1} 
= F + A (\int_{0}^{T} dW_t Y_t')  (\int_{0}^{T} Y_{t} Y_{t}' dt )^{-1} $
when $(\int_{0}^{T} Y_{t} Y_{t}' dt )$ is invertible and, in this case, the 
estimator is unbiased (as the expectation of the second term is zero).
We show here that $\hat{F}_T$ is a consistent and an 
asymptotically efficient estimator of $F$,
 irrespective of the stationarity
 (or stability) of the process, provided $F$ and $A$ together satisfy a
 {\em RANK}
 condition (a), given in Section \ref{sec.2}.  
This RANK condition is essential to prove that 
$ (\int_{0}^{T} Y_{t} Y_{t}' dt ) $ is invertible.
We note here, if $A$ is a nonsingular
 matrix, the {\em RANK} condition automatically holds. 
 In fact, it is also easy to see that for a continuous autoregressive process 
 (i.e., CAR(p)),
 the {\em RANK} condition holds.

We also make another assumption, condition (b). It is the distinctness of
the eigenvalues with positive real parts. However, we point out
that this condition can be relaxed with a condition (b') and also
 that if none of the conditions (b) or (b') hold 
it is still possible
 to proceed with the estimation (see the discussion after Remark \ref{rm2.2}).
 Notice that the condition (b') holds
for the drift $F$ in CAR(p) processes.

The estimation of parameters for the stochastic processes have extensively
studied (see for example, Feigin 
\cite{feigin}, 
Basawa, Feigin and Heyde 
\cite{basfeiheyd},
Basawa and Prakasa Rao 
\cite{key11},
Dietz and Kutoyants
\cite{dtzkuto},
 Kutoyants 
\cite{kuto, kuto1}, 
Barndorff-Nielson and S\/orensen 
\cite{barnsoren},
Kutoyants and Pilibossian 
\cite{kutopili},
Jankunas and Khasminskii
\cite{jnkh},
Khasminskii, Krylov and Moshchuk 
\cite{khkm}
Prakasa Rao
\cite{prakasa1, prakasa2}
 and references therein). 
Therefore,
the estimation of the paramater and its asymptotic studies have not been new.
However, as far as we know, full study of multidimensional OU processes parameter estimation and the study of its asymptotics have not been done for the mixed model.
Apart from 
 showing 
consistency
and asymptotic efficiency for the multidimensional (matrix valued)
variable 
that does not follow from that of univariate or vector valued case
(see, for example, Kaufmann \cite{kauf}, Wei \cite{key15},
Basawa and Prakasa Rao 
\cite{key11},
Dietz and Kutoyants
\cite{dtzkuto},
 Kutoyants 
 \cite{kuto, kuto1}, 
Barndorff-Nielson and S\/orensen 
\cite{barnsoren},
 Prakasa Rao
\cite{prakasa1, prakasa2} and references therein)
it also develops new methodology 
to deal with such cases
as is done in Kaufmann 
\cite{kauf} and Wei \cite{key15}.

Our paper is organized as follows. In Section \ref{sec.2}, we present the
basic assumptions and the main theorems. 
In Section \ref{sec.3}, 
 we describe the case in which the eigenvalues of
$F$ have positive real parts.
Methodology used here is similar to that of Lai and Wei's paper \cite{key1},
 while the case
in which the eigenvalues of $F$ have negative or zero real parts
is quite different from them and it is discussed in Section \ref{sec.4}.
 This case,
in fact, combines the three cases, zero eigenvalues,
 purely imaginary eigenvalues
and the eigenvalues with negative real parts. Details on the rates
of growth and so forth for zero eigenvalues
and imaginary eigenvalues are given in the Appendix.
Section \ref{sec.5}
 examines
 the mixed case 
 for consistency.
The section \ref{sec.6} 
presents
 the results on asymptotic efficiency and
 some concluding remarks.

\section{Basic Assumptions and the Main Theorem}
\setcounter{equation}{0}
\setcounter{section}{2}
\label{sec.2}
\vspace{-0.4cm}
We can decompose
any $p \times p$ matrix $F$ into
the rational canonical form
\begin{eqnarray}
MF = GM = \left( \begin{array}{cccccc} G_0 & 0 \\ 0 & G_1 \end{array}\right)
\left( \begin{array}{cccccc} M_0 \\ M_1 \end{array}\right) \nn
\end{eqnarray}
where $G_i$ are $p_i \times p_i$ matrices and $M_i$ are $p_i \times p$
 matrices for 
$i=0,1$ and $p_0+p_1 = p$. Rows of $M_i$ and rows of $M_j$ are orthogonal for
 $i\not= j$.\\
 All roots of $G_0$ lie in the right half space; all roots of
 $G_1$ lie
on the left half space.

 {\bf EXAMPLE} Let
\begin{eqnarray}
A=\left( \begin{array}{cccccc}
2 & -1 & 0 & 1 & 0 \\
0 & -8 & 6 & 14 & 1 \\
0 & 10 & -4 & -14 & -1 \\
0 & -10 & 6 & 16 & 1 \\
0 & -5 & 3 & 7 & 0 
\end{array}\right) . \nonumber
\end{eqnarray}
Then the characteristic polynomial of A is
\begin{eqnarray}
f(t)=(t-2)^3(t^2+1). \nonumber
\end{eqnarray}
Thus $\phi_1 (t)=t-2$ and $\phi_2 (t)=t^2+1 $ are the distinct irreducible
 monic divisors of $f(t)$. After computation, we find that
  $g(t)=\phi_1(t)^2\phi_2(t)=(t-2)^2(t^2+1)$ is the minimal polynomial of $A$
 and thus the companion matrices for $\phi_1^2(t)=(t-2)^2$ and $\phi_1(t)=t-2$
 are given by
\begin{eqnarray}
\left( \begin{array}{cccccc} 0 & -4 \\ 1 & 4 \end{array}\right) \qquad {\rm  and} \qquad 2 . \nonumber 
\end{eqnarray}
Similarly, the companion matrix for $\phi_2(t)=t^2+1$ is
\begin{eqnarray}
\left( \begin{array}{cccccc} 0 & -1 \\ 1 & 0 . \end{array}\right) \nonumber
\end{eqnarray}
The
 rational canonical form of $A$ is
thus
\begin{eqnarray}
H_A=
\left( \begin{array}{cccccc}
0 & -4 & 0 & 0 & 0 \\
1 & 4 & 0 & 0 & 0 \\
0 & 0 & 2 & 0 & 0 \\
0 & 0 & 0 & 0 & -1 \\
0 & 0 & 0 & 1 & 0 
\end{array}\right) \nonumber
\qed \end{eqnarray}

\vspace{0.3cm}
 In the example above, the rational canonical form of $A$ is
 formed by
 3 blocks: $\left( \begin{array}{cccccc} 0 & -4 \\ 1 & 4\end{array}\right)$,
 $2$ and
  $\left( \begin{array}{cccccc} 0 & -1 \\ 1 & 0 \end{array}\right)$. Therefore
 the dimensions of the
   3 blocks are 2, 1 and 2 respectively.

{\bf ASSUMPTION} 
\begin{eqnarray}
\label{eq1.2}
(a) \qquad 
{\rm RANK} \left(\left[  A: FA: \cdots :F^{p -1}A  \right]\right) = p .
\end{eqnarray}

\medskip
(b) The eigenvalues of $F$, which have positive real parts, are all distinct.

Observe that, from (\ref{eq1.1})
$Y_t = e^{Ft} Y_0 + \int_{0}^{t} e^{F(t-s)} A dW_s$ and thus
have a multivariate Gaussian distribution with the mean $e^{Ft}$ and the
covariance matrix $\int_{0}^{t} e^{Ft} A A' e^{F't}$. 
Since $Y_t$ is Gaussian it has a positive density if and only if the 
covariance matrix is nonsingular. 
The RANK assumption which is the special case of 
H\"{o}rmander's hypoellipticity condition ensures the positive density
of $Y_t$ (for details, see \cite{horm}), and hence the nonsingularity of covariance matrix.

Following Basawa and Rao (\cite{key11}, pp.) it is clear that
$\int_0^T Y_t Y_t' dt$ is nonsingular under the RANK assumption.

Let $F_A =[A:FA: \cdots :F^{p-1}A]$. Then RANK $(F_A)= p$ by the 
RANK assumption.
 Consider
 for $i=0,1$,
\begin{eqnarray}
p_i= {\rm RANK}(M_iF_AF_A^{-1}) \leq {\rm RANK}(M_iF_A) \leq p_i \nn
\end{eqnarray}
where $F_A^{-1}$ is the right inverse of $F_A$.
Therefore, RANK $(M_iF_A)= p_i$ for $i=0,1$. \\
Since
\begin{eqnarray}
M_i F_A &=& \left[ M_i \left[ A : FA : \cdots : F^{p-1}A \right] \right] \nonumber \\
&=& \left[ M_iA : M_iFA : \cdots : M_iF^{p-1}A \right] \nonumber \\
&=& \left[ M_iA : G_iM_iA : \cdots : G_i^{p-1}M_iA \right], \nn
\end{eqnarray}
and as the higher power of $G_i$ can be expressed as a linear combination
of $I$, $G_i, \dots, G_i^{p_i-1}$,  
\begin{eqnarray}
{\rm RANK} \left[ M_iA : G_iM_iA : \cdots : G_i^{p_i-1}M_iA \right]
 = {\rm RANK} \left[ M_iA : G_iM_iA : \cdots : G_i^{p-1}M_iA \right]
  = p_i . \ \ 
\label{eq1.10}
\end{eqnarray}
If we transform the process $Y_t$ to $U_{it}=M_iY_t$ for $i=0,1$,
\begin{eqnarray}
M_idY_t &=& M_iF Y_t dt + M_iAdW_t , \nn\\
i.e., \ \ dU_{it} &=& G_i U_{it} dt + \left( M_i A \right) dW_t . \nn
\end{eqnarray}
From (\ref{eq1.10}) and the argument given above, we conclude that
$\int_0^T U_{it}U_{it}' dt$
is positive definite a.s. for $i=0,1$.\\

We now present our main theorems whose proofs are given in Section \ref{sec.5} and
in Section \ref{sec.6}, respectively.
Throughout the paper, we use $\la_{\min} (C)$ and $\la_{\max} (C)$
to denote
 the minimum and maximum  eigenvalues of a matrix $C$. 
\begin{theorem}
\label{th1.1}
 Suppose, for the Ornstein-Uhlenbeck process defined in
(\ref{eq1.1}), the
 assumptions (a) and (b)
 hold.
 Define $\hat F_T = (\int_0^T dY_t Y_t')(\int_0^T Y_t Y_t' dt )^{-1}$.
 Then
\begin{eqnarray}
\label{eq1.13}
{\lim \inf}_{T \rightarrow \infty} \frac{1}{T} \lambda_{\min}
 \left( \int_0^T Y_t Y_t' dt \right) >0 \qquad {\rm  a.s. }
\end{eqnarray}
and
\begin{eqnarray}
{\lim}_{T \rightarrow \infty} \hat F_T = F \qquad {\rm  a.s.} \nn
\end{eqnarray}
\end{theorem}

\begin{theorem}
\label{th1.2}
Under the assumptions of Theorem \ref{th1.1}, it follows that \\
$E({\rm Tr}[(\hat{F}_{T} - F) E(C_{T}) (\hat{F}_{T} - F)'])^{1/2} = O(1)$
as $T \to \infty$, where $\hat{F}_{T}$ is as defined in Theorem \ref{th1.1} 
and $C_{T} = \left(\int_{0}^{T} Y_t Y_t' dt \right)$.
\end{theorem}

\section{Eigenvalues in the Right Half Space}
\setcounter{equation}{0}
\setcounter{section}{3}
\label{sec.3}
\vspace{-0.4cm}
We consider the case where all the eigenvalues of $F$ have positive real parts.
In this case, it can be seen that $\|Y_t\| \to \infty$ exponentially fast
as $t \to \infty$. To introduce the main result of this section we
define a Gaussian random variable
$$ Z = Y_0 + \int_0^{\infty} e^{-Fs} A dW_s .
$$
Since all the eigenvalues of $F$ have positive real parts,
it is clear that, $ \ e^{-Ft} Y_t = Y_0 + \int_0^t e^{-Fs} A dW_s \ $ 
converges a.s. to $Z$ as $t \to \infty$. We now derive the following results.

\begin{theorem}
\label{th2.1}
In addition to the assumptions and notations of Theorem \ref{th1.1}, assume
further that
real parts of all the eigenvalues of $F$ are
positive. Then,
\begin{eqnarray}
e^{-FT} \left( \int_0^T Y_t Y_t' dt \right) e^{-F'T} \ \ \
{\rm  converges \ a.s. \ to } \ \ \
 B = \int_0^\infty e^{-Ft} (ZZ') e^{-F't} dt. \nn
\end{eqnarray}
Moreover, $B$ is positive definite with probability 1. Consequently,
\begin{eqnarray}
\label{eq2.2}
{\lim}_{T \rightarrow \infty} T^{-1} \log \lambda_{\min}
 \left( \int_0^T Y_t Y_t' dt
\right) = 2 \lambda_0 \qquad  {\rm  a.s.} \nn\\
{\lim}_{T \rightarrow \infty} T^{-1} \log \lambda_{\max}
 \left( \int_0^T Y_t Y_t' dt
\right) = 2 \Lambda_0 \qquad {\rm  a.s.}
\end{eqnarray}
\end{theorem}

Here and throughout the paper, $\log x$ means the natural logarithm of $x$.
Also, in the sequel
 we shall let $|| x ||$ denote the Euclidean norm of a 
$p$-dimensional
 vector $x=(x_1, \cdots ,x_p)'$, i.e., $||x||^2 =x'x$. Moreover, by viewing a
  $p \times p$ matrix $A_0$ as linear operator, we define
$||A_0|| = \sup_{||x||=1} ||A_0 x||$. Thus, $||A_0||^2$ is equal to the maximum
  eigenvalue of $A_0'A_0$. Moreover, if $A_0$ is symmetric and non-negative
 definite,
 then $||A_0||=\lambda_{\max}(A_0)$. In particular, for the companion matrix
      $e^{-FT}$ in Theorem \ref{th2.1}, we have the following Lemma.

\begin{lemma}
\label{lm2.1}
Under the hypothesis of Theorem \ref{th2.1}
\begin{eqnarray}
\label{eq2.3}
\log ||e^{FT}|| & \sim & \log ||e^{F'T}|| \sim \Lambda _0T,
\nn\\
 {\rm and} \qquad \log ||e^{-FT}|| &\sim& \log ||e^{-F'T}|| \sim -\lambda_0T
\end{eqnarray}
where we use the notation $f(T) \sim CT^k$ to denote
 ${\lim}_{T \rightarrow \infty} T^{-k} f(T) =C$. 
\end{lemma}

 {\bf Proof.}
 Suppose Re$[\lambda_k(F)]>0$ for $k=1,2,\cdots, p$. Then
\begin{eqnarray}
|e^{\lambda_k(F)}|=e^{{\rm Re}[\lambda_k(F)]} >1 \qquad {\rm  for } \ \ 
 k=1,2, 
\cdots ,p. \nn
\end{eqnarray}
Let $\lambda_0 =\min_{1 \leq k \leq p} {\rm  Re } [\lambda_k (F)]$,
 $\Lambda_0=\max_{1 \leq k 
\leq p} {\rm  Re } [\lambda_k (F)]$. Denote 
the spectral radius of
 $F$ by $r_\sigma(F)$ (cf. \cite{key13}). Then
\begin{eqnarray}
{\lim}_{T \rightarrow \infty} ||e^{FT}||^{\frac{1}{T}} = r_\sigma (F)
= {\sup}_{\lambda \in \sigma (e^F)} | \lambda | 
= \exp \left[ {\sup}_{\lambda \in \sigma (F)} {\rm  Re } (\lambda ) \right] 
= e^{\Lambda_0} \nn
\end{eqnarray}
and so $\log || e^{FT}|| \sim \log ||e^{F'T}|| \sim \Lambda_0T$.
Similarly, $\log ||e^{-FT}|| \sim \log ||e^{-F'T}|| \sim -\lambda_0 T$ since
\begin{eqnarray}
{\lim}_{T \rightarrow \infty} ||e^{-FT}||^{\frac{1}{T}} &=& 
{\sup}_{ \lambda \in \sigma (e^{-F})} | \lambda | 
= \exp \left[ {\sup}_{\lambda \in \sigma ( -F)} {\rm  Re } (\lambda ) \right]
= e^{-\lambda_0}. \nn
\end{eqnarray}
Thus, we have the proof of Lemma \ref{lm2.1} \qed

\vspace{0.2cm}
{\bf Proof of Theorem \ref{th2.1}.} 
Let $Z_t =Y_0 + \int_0^t e^{-Fs}AdW_s$, then $Y_t =e^{Ft}Z_t$ and
\begin{eqnarray}
Z_t \ \ {\rm converges \ a.s. \ to } \ \  Z=Y_0 + \int_0^\infty e^{-Fs}AdW_s.
\nn
\end{eqnarray}
Let $B_T=\int_0^T e^{-Ft}Z_TZ_T'e^{-F't}dt$,
\begin{eqnarray}
\label{eq2.4}
&& \left| \left| e^{-FT} \left( \int_0^T Y_t Y_t' dt \right) e^{-F'T} -B_T \right| \right| \nonumber \\
&=& \left| \left| \int_0^T e^{-F(T-t)}Z_t Z_t'e^{-F(T-t)}dt -\int_0^T e^{-Ft}Z_TZ_T'e^{-F't}dt \ \right| \right| \nonumber \\
&=& \left| \left|  \int_0^T e^{-Ft} \left(Z_{T-t}Z_{T-t}'-Z_TZ_T'\right) e^{-F't}dt \right| \right| \nonumber \\
&\leq& \int_0^T ||e^{-Ft}|| \hspace{0.1cm} ||e^{-F't}|| \left( || Z_{T-t} || + ||Z_T|| \right) ||Z_T-Z_{T-t}||dt  \nonumber \\
&=& \int_0^{T/2}  ||e^{-Ft}||^2
 (||Z_{T-t}||+||Z_T||) ||Z_T-Z_{T-t}|| dt \nonumber \\
&& + \int_{T/2}^T ||e^{-Ft}||^2
 ( ||Z_{T-t}|| + ||Z_T ||) ||Z_T -Z_{T-t}|| dt .
\end{eqnarray}
Since $Z_t$ converges almost surely to a finite random variable $Z$,
$\sup_{\{t \ge 0\}} \|Z_t\|$ is finite almost surely
 and for each $t \ge T/2$, 
 $||Z_T -Z_{T-t}||$, being a cauchy sequence, converges to zero,
 almost surely, as $T \to \infty$. Also, 
 by Lemma \ref{lm2.1},
 $ \int_{0}^{\infty}  ||e^{-Ft}||^2 dt < \infty$.
Thus, we get,
 $\forall \omega$ outside a null set, $\forall \ep > 0$,
there exists a $T_{0}(\omega)$ such that
$\| Z_t(\omega) - Z(\omega) \| < 
\ep/(1 + \int_{0}^{\infty}  ||e^{-Ft}||^2 dt + 2 \sup_{\{t \ge 0\}} \|Z_t
 (\omega)\|)$
for all $t \geq T_{0}(\omega)$.
Fixing one such $\omega$, for $T \geq 2 T_{0}(\omega)$
we have the first integral of (\ref{eq2.4}), which is less
than $\ep$ and the second integral goes to zero as
$\sup_{\{t \ge 0\}} \|Z_t (\omega)\|$ is finite and 
$\int_{T/2}^{T}  ||e^{-Ft}||^2 dt \to 0$ as $T \to \infty$.

Let $B=\int_0^\infty e^{-Ft} ZZ'e^{-F't}dt$, then with probability 1,
\begin{eqnarray}
\label{eq2.5}
&& ||B_T -B || \nonumber \\
&\leq & \int_T^\infty ||e^{-Ft}ZZ'e^{-F't} ||dt + \int_0^T ||e^{-Ft}(ZZ'-Z_TZ_T' ) e^{-F't}||dt \nonumber \\
&\leq& ||ZZ'|| \int_T^\infty ||e^{-F't}|| \hspace{0.1cm} ||e^{-Ft}||dt + ||ZZ'-Z_TZ_T'||\int_0^T ||e^{-Ft}||\hspace{1mm}||e^{-F't}|| dt  \nonumber \\
&\rightarrow& 0 {\rm  \quad a.s., \quad as \quad } T \rightarrow \infty .
\end{eqnarray}
Therefore,
\begin{equation}
\label{eq2.6}
e^{-FT} \left( \int_0^T Y_t Y_t'dt \right) e^{-F'T} 
\ \ {\rm converges \ a.s. \ to } \ \  B=\int_0^\infty e^{-Ft} ZZ' e^{-F't}dt.
\end{equation}

To show $B= \int_0^\infty e^{-Ft} ZZ' e^{-F't} dt$ is positive definite with
 probability 1, 
observe that $Z$ has positive Gaussian density. Hence $P(Z \neq 0) = 1$.
Fix an $\omega$, such that $Z(\omega) \neq 0$.
Suppose, if possible,
\begin{eqnarray}
x'\left( \int_0^\infty e^{-Ft}Z(\omega)Z(\omega)'e^{-F't} dt \right) x =0 \quad {\rm  for \
 some \ nonzero \ vector } \ \ x \in {\cal R}^p  
 .\nn
\end{eqnarray} 
Then, for almost all $t \in (0,T)$, $x'e^{-Ft}Z(\omega)=0$, 
i.e., for almost all $t \in (0,T)$,
$\sum\limits_{k=0}^\infty \frac{1}{k!} (-1)^k x'F^k t^k Z(\omega) =0$.
This implies
 $x'F^kZ(\omega)=0$, 
 for $k = 0, 1, \cdots ,p-1$.
By the assumption (b), $\sum_{k=0}^{p-1} a_k F^k$ is nonsingular for
any real number $a_k$ with not all of them being zero. Hence, for
any nonzero vector in ${\cal R}^p$, in particular for $x$,
 $\ x'\sum_{k=0}^{p-1} a_k F^k$ is 
a nonzero vector. In other words, for nonzero vector $x$,
  $\sum_{k=0}^{p-1} a_k (x' F^k)$ is nonzero for any nonzero vector
$(a_0, \ldots, a_{p-1})$. Thus
$\left( \begin{array}{cccccc} x' \\ x'F \\ \vdots \\ x'F^{p-1} 
\end{array}\right) $ is a nonsingular matrix.
Hence,
$ \ \  \left( \begin{array}{cccccc} x' \\ x'F \\ \vdots \\ x'F^{p-1} 
\end{array}\right) Z(\omega) = 0$ 
\  implies \
 $Z(\omega) = 0$,
which is a contradiction.
Thus, we arrive at a contradiction since $Z$ has a positive Gaussian density
and hence $Z$ cannot be equal to zero on a set of positive measures.
Therefore, we conclude that
$B$ is positive definite with probability one.

To prove (\ref{eq2.2}), we state the following elementary results
(for the proof, see Lemma 2 of \cite{key1}):

\begin{lemma}
\label{lm2.2}
 Let $A$, $C$ be $p \times p$ matrices such that $C$ is symmetric and
 non-negative definite. Then
\begin{eqnarray}
\lambda_{\max}(C)\lambda_{\max}(AA') \geq \lambda_{\max}(ACA') \geq \lambda_{\min}(C)\lambda_
{\max}(AA'), \nn\\
\lambda_{\max} (C) \lambda_{\min}(AA') \geq \lambda_{\min}(ACA') \geq \lambda_{\min} (C) \lambda_{\min}(AA'). \nn
\end{eqnarray}.
\end{lemma}

We continue the proof of (\ref{eq2.2}) of Theorem \ref{th2.1}.
  From Lemma \ref{lm2.2} we get,
\begin{eqnarray}
 \log \lambda_{\min} \left( \int_0^T Y_t Y_t' dt \right) &\leq & \log \lambda_{\max} \left[ e^{-FT} \left( \int_0^T Y_t Y_t' dt \right) e^{-F'T} \right]  \nonumber \\
&& - \log \lambda_{\max} \left( e^{-FT} e^{-F'T} \right)  \nonumber \\
& \sim & 2\lambda_0T. \nn
\end{eqnarray}
Also,
\begin{eqnarray}
\log \lambda_{\min} \left( \int_0^T Y_t Y_t' dt \right) & \geq & \log \lambda_{min} \left[ e^{-FT} \left( \int_0^T Y_t Y_t' dt \right) e^{-F'T} \right] \nonumber \\
&& + \log \lambda_{\min} \left( e^{FT} e^{F'T} \right) \nonumber \\
& \sim & 2\lambda_0T. \nn
\end{eqnarray}
Therefore
\begin{eqnarray}
{\lim}_{T \rightarrow \infty} \frac{1}{T} \log \lambda_{\min} \left( \int_0^T Y_t Y_t
' dt \right) = 2\lambda_0 \qquad {\rm  a.s.} \nn
\end{eqnarray}

On the other hand,
\begin{eqnarray}
\log \lambda_{\max} \left( \int_0^T Y_t Y_t' dt \right)
& \leq & \log \lambda_{\max} \left[ e^{-FT} \left (\int_0^T Y_t Y_t' dt\right) e^{-F'T} \right] \nonumber \\
&&  + \log \lambda_{\max} \left( e^{FT}e^{F'T} \right) \nonumber \\
&\sim & 2\Lambda_0T. \nn
\end{eqnarray}
Also,
\begin{eqnarray}
\log \lambda_{\max} \left( \int_0^T Y_t Y_t' dt \right) 
& \geq & \log \lambda_{\min} \left[ e^{-FT} \left( \int_0^T Y_t Y_t' dt \right) e^{-F'T} \right] \nonumber \\
&& - \log \lambda_{\min} \left( e^{-FT}e^{-F'T} \right) \nonumber \\
& \sim &2\Lambda_0T . \nn
\end{eqnarray}
Therefore
\begin{eqnarray}
{\lim}_{T \rightarrow \infty} \frac{1}{T} \log \lambda_{\max} \left( \int_0^T Y_t Y_t
' dt \right) = 2\Lambda_0 \qquad {\rm a.s.} \nn
\end{eqnarray}       
Hence, we have the proof of Theorem \ref{th2.1}. \qed 

\begin{corollary}
\label{cor2.1}
 Under the same assumptions and notations as in Theorem \ref{th2.1},
\begin{eqnarray}
(i) && {\lim}_{T \rightarrow \infty} \int_0^T ||e^{-FT}Y_t ||dt =
 \int_0^\infty ||e^{-Ft}Z || dt < \infty \qquad {\rm  a.s.} \\
(ii) && \frac{1}{\sqrt{T}}\left( \int_0^T dW_t Y_t' \right) e^{-F'T} =
 O(T^{-1/2}) . \nn
\end{eqnarray}
\end{corollary}

{\bf Proof.} (i) \hspace{4pt}
 Given $\epsilon >0, \forall \omega$ outside a
 null set, $\exists T_{0}(\omega )$ such that
\begin{eqnarray}
|| Z_t - Z || < \epsilon \qquad \forall t \geq T_{0}( \omega ). \nn
\end{eqnarray}
For $T> T_{0}(\omega)$,
\begin{eqnarray}
&& \left| \int_0^T ||e^{-F(T-t)}Z_t ||dt - \int_0^T ||e^{-F(T-t)}Z||dt \right|
 \nonumber \\
& \leq & \int_0^T ||e^{-F(T-t)}Z_t -e^{-F(T-t)}Z ||dt \nonumber \\
& \leq & \int_0^T ||e^{-F(T-t)}|| \hspace{0.1cm} ||Z_t -Z|| dt \nonumber \\
& \leq & \int_0^{T_{0}(\omega )}||e^{-F(T-t)}|| \hspace{0.1cm} ||Z_t - Z|| dt + 
\int_{T_{0}(\omega )}^T ||e^{-F(T-t)}||\hspace{0.1cm} ||Z_t -Z|| dt . \nn
\end{eqnarray}
As $T \rightarrow \infty$, the first term tends to 0 since $||e^{-F(T-t)}||
 \rightarrow 0$. The second term also tends to 0 since $Z_t \rightarrow Z$ and
 $\int_{T_{0}(\omega )}^T ||e^{-F(T-t)}||dt
 \le \int_0^T ||e^{-F(T-t)}||dt =  \int_0^T ||e^{-Ft}||dt \le  \int_0^{\infty}
 ||e^{-F(T-t)}||dt $,
which is finite. Therefore,
\begin{eqnarray}
{\lim}_{T \rightarrow \infty} \int_0^T ||e^{-FT}Y_t ||dt &=&
 {\lim}_{T \rightarrow \infty} \int_0^T ||e^{-F(T-t)}Z_t ||dt \nonumber \\
&=& {\lim}_{T \rightarrow \infty} \int_0^T ||e^{-F(T-t)}Z|| dt \nonumber \\
&=& \int_0^\infty ||e^{-Ft}Z|| dt , \nn
\end{eqnarray}
which is finite
 almost surely,
 by Lemma \ref{lm2.1}.

(ii) Let $M_t =\left( \int_0^t dW_s Y_s' \right) e^{-F'T}$,
 which is a square integrable martingale for $0 \leq t \leq T$, with quadratic
 variation,
\begin{eqnarray}
<M>_t=e^{-FT}\left( \int_0^t Y_s Y_s' ds \right) e^{-F'T} =e^{-FT} C_t e^{-F'T}
\nn
\end{eqnarray}
where $C_t = \int_0^t Y_s Y_s' ds$.
By Karatzas and Shreve (cf \cite{key3} p174),
\begin{eqnarray}
\left( \int_0^t dW_s Y_s' \right) e^{-F'T} &=& M_t = B_{<M>_t} \nn\\
&=& O \left( \lambda_{\max} \left( e^{-FT}C_te^{-F'T} \right) \sqrt{ \ln \ln \lambda_{\max} \left( e^{-FT} C_t e^{-F'T}\right)} \right) \nonumber \\
&=& O(1) \nn
\end{eqnarray}
since for $t \leq T, \ \ ||e^{-FT}C_te^{-F'T}|| \leq ||e^{-FT}C_Te^{-F'T}||
 \rightarrow B$, almost surely, as $T \to \infty$
 and $B = O(1)$. 
Therefore,
\begin{eqnarray}
\frac{1}{\sqrt{T}} \left( \int_0^T dW_t Y_t' \right) e^{-F'T} = O(T^{-1/2}) \nn
\end{eqnarray}
This completes the proof of Corollary \ref{cor2.1}.
 \qed 

\begin{remark}
\label{rm2.1}
\end{remark}
If all the eigenvalues of $F$ have positive real parts, we can relax 
condition (b) by
\begin{eqnarray}
\label{eq1.2a}
 & (b') & \quad \sum_{k=0}^{p-1} a_k F^k \ \ {\rm being \ nonsingular} 
 {\rm \ for \ any \ reals \ \ } a_1, \ldots, a_n {\rm \ \ with \ 
at \ least } \nn\\
&& {\rm  \ one \ of \ them \ being \ nonzero.} 
\end{eqnarray}
Notice that (b') could hold even if all the eigenvalues of $F$ are equal
(say, $\la_0$), but
 the degree of the minimal polynomial of $F$ and the degree of 
the characteristic polynomial of $F$ are equal.

\begin{remark}
\label{rm2.2}
Suppose, assumption (b) does not hold. One can still estimate the eigenvalues of $F$.
\end{remark}

Let the characteristic polynomial of $F$ be given as
$\phi_{F}(x) = a_0\Pi_{i=1}^k (x - \la_i)^{p_i}
\Pi_{j=1}^l (x^2 + b_j x + c_j)^{q_j}$
 where $\la_i$ are the real roots of
 multiplicity $p_i$ and $x^2 + b_j x + c_j$ are the irreducible polynomials
giving the complex roots with multiplicity $q_j$ and $a_0$ is a constant.
Let the minimal polynomial
of $F$ be given by
$\psi_F(x) = 
\Pi_{i=1}^k (x - \la_i)^{r_i}
\Pi_{j=1}^l (x^2 + b_j x + c_j)^{s_j}$ 
with $r_i \le p_i$ and $s_j \le q_j$.
If $r_i = p_i$ and $s_j = q_j$ for all $i, \ j$, then
the degree of the minimal polynomial of $F$ and the 
degree of the characteristic polynomial of $F$ are the same and 
the assumption (b') holds and our results follow.
 If some of the $r_i$s are less than $p_i$s
and/or $s_j$s are less than $q_j$,  
then, (b') does not hold for $F$.
However, in that case, one can transform $F$ in the rational canonoical
form as 
\begin{eqnarray}
\left( \begin{array}{cccccc} J_1 \\ \vdots \\ J_k \\ K_1 \\ \vdots \\ K_l \\ L \end{array}\right) F 
=
 \left( \begin{array}{cccccccc}
  B_1 & \cdots & 0 & 0 & \cdots & 0 & 0 \\
   0 & \ddots & 0 & 0 & \cdots & 0 & 0 \\
   0 & \cdots & B_k & 0 & \cdots & 0 & 0 \\
   0 & \cdots & 0 & C_1 & \cdots & 0 & 0 \\ 
   0 & \cdots & 0 & 0 & \ddots & 0 & 0 \\ 
   0 & \cdots & 0 & 0 & \cdots & C_l & 0 \\ 
   0 & \cdots & 0 & 0 & \cdots &  0 & D 
  \end{array}\right)
\left( \begin{array}{cccccc} J_1 \\ \vdots \\ J_k \\ K_1 \\ \vdots \\ K_l \\ L \end{array}\right) 
=
\left( \begin{array}{cccccc} B_1 J_1 \\ \vdots \\ B_k J_k \\ C_1 K_1 \\ \vdots \\ C_l K_l   \\ D L \end{array}\right) \nn
\end{eqnarray}
where $J_i$, $K_j$ and $L$ are rectangular matrices of full row rank,
$(p_i-r_i)$, $(q_j - s_j)$, $(\sum_i r_i + \sum_j s_j)$, respectively,
and $D$ is 
a square matrix of the dimension the same as the degree of the
 minimal polynomial of $F$ $({\rm i.e., \ same \ as \ }
(\sum_i r_i + \sum_j s_j))$.
For each $j$, $C_j$
is a partitioned diagonal matrix
(i.e., only the diagonal blocks are nonzero blocks),
each block is of dimension $2 \times 2$, and its diagonal block matrices are
identical 
and repeating exactly $(q_j - s_j)$ times and
have the characteristic polynomial
$x^2 + b_j x + c_j$,
and, for each $i$, $B_i$ is a diagonal matrix 
with diagonal entries consisting of the real characteristic root $\la_i$
 repeating exactly $(p_i - r_i)$ times.
 Thus, we can work with $D$ instead of $F$.
For $D$ the 
assumption (b') holds, since the degree of minimal polynomial of $D$ is same
as that of $F$ and, consequently,
 the degree of the minimal polynomial of $D$ is the same as the
degree of the characteristic polynomial of $D$.
Estimation of $D$ can be done using the SDE of $LY_t$.
For $B_i$ and $C_j$, one can consider each one 
 separately
and transform $Y_t$ to $J_i Y_t$ and $K_j Y_t$ and use the SDE of any
component of $J_i Y_t$ (as it has the Markov property) to estimate $\la_i$ 
and the SDE of the
 first two (or, any (2m-1)th and 2mth)
 components of $K_j Y_t$ together,
 as they have the Markov property,
 to estimate a diagonal block of $C_j$. Hence the assertion in the last remark.

\section{Eigenvalues on the Left Half Space}
\setcounter{equation}{0}
\setcounter{section}{4}
\label{sec.4}
\vspace{-0.4cm}
 In this Section, we study the asymptotic behavior of
OU processes where the real parts of all the eigenvalues of $F$ are either
zero or negative. Unlike the exponential rate of growth for $||Y_T||$,
 $\lambda_{\max} ( \int_0^T Y_t Y_t' dt )$, $\lambda_{\min}
 (\int_0^T Y_t Y_t' dt )$ in Theorem \ref{th2.1} and Corollary \ref{cor2.1}
 for the the process where all the eigenvalues of $F$ have positive real parts,
the following theorem shows that these quantities grow at
 most polynomially fast in $t$ for these processes.

 For stable processes $Y_t$ (i.e., eigenvalues of $F$ 
with negative real parts), we know from Basak and Bhattacharya
 \cite{key10} that
\begin{eqnarray}
| Y_t^x - Y_t^0 | \rightarrow 0
 {\rm  \quad a.s. \quad as \quad } t \rightarrow \infty. \nn
\end{eqnarray}
Therefore, the property of $Y_t$ starting at $x$ is the same as that from $0$.
Hence, without loss of generality, we can assume that $Y_0=0$. \\

\begin{theorem}
\label{th3.1}
 Suppose,
for the Ornstein-Uhlenbeck process defined in
(\ref{eq1.1}), the RANK condition (\ref{eq1.2}) holds
and all the eigenvalues of $F$ have negative real parts.
 Then  
\begin{eqnarray}
\label{eq3.2}
 && {\lim \inf}_{T \rightarrow \infty} \frac{1}{T} \lambda_{\min}
 \left( \int_0^T Y_t Y_t' dt \right) > 0 \qquad {\rm  a.s.}
\end{eqnarray}
Moreover,
\begin{eqnarray}
\label{eq3.3}
\lambda_{\max} \left( \int_0^T Y_t Y_t' dt \right) &=& O(T) \hspace{1cm} {\rm  a.s.}  
\end{eqnarray}
\end{theorem}

{\bf Proof.} 
To prove (\ref{eq3.2}) and (\ref{eq3.3}), consider each component
 $Y_t^i, Y_t^j$ of $Y_t$, $i,j =1, \cdots ,p$. Let $\pi$ be the invariant
 distribution of $Y$. Then by the Strong Law of Large Numbers,
\begin{eqnarray}
\frac{1}{T} \int_0^T Y_t^i Y_t^j dt \rightarrow E_\pi (Y^iY^j) < \infty
 {\rm  \quad as } \ \ T \rightarrow \infty , \nn
\end{eqnarray}
which follows, afortiori, by the Law of the Iterated Logarithm by
 Basak \cite{key9}. 
Therefore,
\begin{eqnarray}
\frac{1}{T} \int_0^T Y_t Y_t' dt \rightarrow
 E_\pi (Y Y') =\int_0^\infty e^{Fu}AA'e^{F'u} du , \nn
\end{eqnarray}
which is positive definite a.s. Therefore,
\begin{eqnarray}
&& {\liminf}_{T \rightarrow \infty} \frac{1}{T} \lambda_{\min}
 \left( \int_0^T Y_t Y_t' dt \right) >0 \qquad {\rm  a.s.} \nonumber  \\
{\rm  and } \qquad && \lambda_{\max} \left( \int_0^T Y_t Y_t'dt \right) = O(T)
 \qquad \qquad {\rm  a.s.} \nn
\end{eqnarray}
Hence, the proof. \qed 

\begin{remark}
\label{rm3.1}
\end{remark}
 (i) It is not difficult to see that for stable $Y_t$, for any $m \geq 1$,
 $E \left[ {\sup}_{k-1 \leq t \leq k}(Y_t'PY_t)^m\right]$ is bounded
 uniformly over k. Hence, it would follow, for any $\delta > 0$,
 $||Y_t||=o(t^{\frac{1}{2m}+\delta})$ a.s. \\
(ii) On the other hand, since $Y_t \rightarrow Y$ in distribution
and $Y$ is finite with probability one,
 one obtains
 $Y_t=O_p(1)$.

\begin{corollary}
\label{cor3.1}
 With the same notations and assumptions as in Theorem 
\ref{th3.1}, let $C_T = \int_0^T Y_t Y_t' dt$. Then 
\begin{eqnarray}
& (i) & \hspace{0.8cm}||C_T^{-1/2}|| = O(T^{-1/2}), \ \ \ {\rm a.s.} \hspace{4cm} \nonumber \\
& (ii) & \hspace{0.8cm} {\lim}_{T \rightarrow \infty}
 Y_T' C_T^{-1} Y_T =0 \qquad {\rm  a.s.} \hspace{4cm} \nn
\end{eqnarray}
\end{corollary}

{\bf Proof.}
 (i) Since ${\lim \inf}_{T \rightarrow \infty}
 \frac{1}{T} \lambda_{\min} (C_T) >0$ a.s. from (\ref{eq3.2}), therefore
\begin{eqnarray}
||C_T^{-1/2}||^2 =\lambda_{\max} (C_T^{-1}) =\frac{1}{\lambda_{\min}(C_T)}
 =O(T^{-1}) \quad {\rm a.s.} \nn
\end{eqnarray}
(ii) 
By the previous remark \ref{rm3.1} (i), we note that,
\begin{eqnarray}
 ||Y_T'C_T^{-1}Y_T|| & \leq & ||Y_T||^2 ||C_T^{-1}|| \nn\\
&=& o(T^{1/2+2\delta})O(T^{-1}) \quad {\rm a.s., \ for \ some \ } \de > 0,
 {\rm \ \ small }
\nn\\
&=& O(T^{-1/2+2\delta})
\nn
\end{eqnarray}
Hence, the proof.  \qed

\begin{theorem}
\label{th3.2}
 Suppose eigenvalues of $F$ have
either negative or zero real parts
(i.e., the eigenvalues are on the Left Half Space, which includes
 zero eigenvalues,
purely imaginary eigenvalues, eigenvalues with negative real parts).
Then,
\begin{eqnarray}
{\lim}_{T \rightarrow \infty} Y_T' 
\left( \int_0^T Y_tY_t' dt \right)^{-1} Y_T =0 \quad {\rm a.s.} \nn
\end{eqnarray}
\end{theorem}

\medskip
To prove Theorem \ref{th3.2}, we need the following lemma:

\begin{lemma}
\label{lm3.1}
 Let $\epsilon > 0$; define $F^\epsilon = F - \epsilon I$ and 
$dY_t^\epsilon =F^\epsilon Y_t^\epsilon dt + AdW_t$.
Then
$\frac{\partial}{\partial \epsilon} \ln \left[ (Y_T^\epsilon)'(C_T^\epsilon)^{-1}(Y_T^\epsilon) \right] $
is bounded below, almost surely,
 uniformly for large values of $T$. 
\end{lemma}

{\bf Proof.} 
Let $\dot Y_t^\epsilon = \frac{\partial}{\partial \epsilon} Y_t^\epsilon$.
 Then we have
\begin{eqnarray}
d \dot Y_t^\epsilon &=& \left( - Y_t^\epsilon + F^\epsilon \dot Y_t^\epsilon \right) dt, \nn
\end{eqnarray}
or jointly,
\begin{eqnarray}
d \left( \begin{array}{cccccc} Y_t^\epsilon \\ \dot Y_t^\epsilon \end{array}\right) &=&
\left( \begin{array}{cccccc} F^\epsilon & 0 \\ -I & F^\epsilon \end{array}\right)
\left( \begin{array}{cccccc} Y_t^\epsilon \\ \dot Y_t^\epsilon \end{array}\right)
dt + \left( \begin{array}{cccccc} A \\ 0 \end{array}\right) dW_t . \nn
\end{eqnarray}
Since all eigenvalues of
 $\left( \begin{array}{cccccc} F^{\epsilon} & 0  \\ -I & F^{\epsilon}
 \end{array}\right)$ have negative real parts,
 $\left( \begin{array}{cccccc} Y_t^\epsilon \\
 \dot Y_t^\epsilon \end{array}\right)$ is stable. Therefore,
\begin{eqnarray}
\left| \left| \left( \begin{array}{cccccc} Y_t^\epsilon \\ \dot Y_t^\epsilon
 \end{array}\right) \right| \right| = 
o(t^{\frac{1}{4} + \delta} ) \quad {\rm a.s. \ for \ some \quad } \de > 0 \nn
\end{eqnarray}
 and
\begin{eqnarray}
\frac{1}{T} \int_0^T \left( \begin{array}{cccccc} Y_t^\epsilon \\ \dot
 Y_t^\epsilon \end{array}\right) \left( \begin{array}{cccccc} Y_t^\epsilon &
 \dot Y_t^\epsilon \end{array}\right) dt \nn
\end{eqnarray}
is positive definite (since the RANK condition holds here as well)
and it converges almost surely to some positive definite constant matrix
 as $T \rightarrow \infty$. Therefore,
 $(C_T^\epsilon)$ and $(\dot C_T^\epsilon)$
 have the same order
 where $C_T^\epsilon=\int_0^T Y_t^\epsilon Y_t^\epsilon  dt$ and
 $\dot C_T^\epsilon =\int_0^T \dot Y_t^\epsilon \dot Y_t^\epsilon dt$. Hence
\begin{eqnarray}
\label{eq3.4}
(\dot C_T^\epsilon)(C_T^\epsilon)^{-1} =O(1) \quad {\rm a.s. \ as \quad }
T \to \infty .
\end{eqnarray}

By Corollary \ref{cor3.1},
\begin{eqnarray}
{\lim}_{T \rightarrow \infty} (Y_T^\epsilon)'(C_T^\epsilon)^{-1}(Y_T^\epsilon)
 &=& 0 \quad {\rm a.s. \ and  } \nonumber \\
{\lim}_{T \rightarrow \infty} (\dot Y_T^\epsilon)'(C_T^\epsilon)^{-1}(\dot Y_T^\epsilon) &=&
 {\lim}_{T \rightarrow \infty} (\dot Y_T^\epsilon)'(\dot C_T^\epsilon)^{-1}(\dot Y_T^\epsilon) =0 \quad {\rm a.s.} \nn
\end{eqnarray}

Consider
\begin{eqnarray}
&& \frac{\partial}{\partial \epsilon}(Y_T^\epsilon)'(C_T^\epsilon)^{-1}(Y_T^\epsilon) \nonumber \\
&=& 2 (\dot Y_T^\epsilon)'(C_T^\epsilon)^{-1} Y_T^\epsilon + (Y_T^\epsilon)'\frac{\partial}{\partial \epsilon}(C_T^\epsilon)^{-1} Y_T^\epsilon \nonumber \\
&=& 2(\dot Y_T^\epsilon)'(C_T^\epsilon)^{-1}Y_T^\epsilon - (Y_T^\epsilon)'(C_T^\epsilon)^{-1} \left[ \frac{\partial}{\partial \epsilon} C_T^\epsilon \right] (C_T^\epsilon)^{-1} Y_T^\epsilon  \nonumber \\
& \geq & -2 \left[ (\dot Y_T^\epsilon)'(C_T^\epsilon)^{-1} (\dot Y_T^\epsilon )\right]^{1/2} \left[ (Y_T^\epsilon)'(C_T^\epsilon)^{-1} (Y_T^\epsilon)\right]^{1/2} \nonumber \\
&& -(Y_T^\epsilon)'(C_T^\epsilon)^{-1} \left[ \int_0^T (Y_u^\epsilon)(\dot Y_u^\epsilon)'du + \int_0^T (\dot Y_u^\epsilon)(Y_u^\epsilon)'du \right] (C_T^\epsilon)^{-1}Y_T^\epsilon \nonumber \\
& \geq & -2 \left[ (\dot Y_T^\epsilon)'(C_T^\epsilon)^{-1} (\dot Y_T^\epsilon) \right]^{1/2} \left[ (Y_T^\epsilon)'(C_T^\epsilon)^{-1}(Y_T^\epsilon)\right]^{1/2}  \nonumber \\
&& -2 \left[ (Y_T^\epsilon)'(C_T^\epsilon)^{-1}(Y_T^\epsilon)\right] \int_0^T \left[ (Y_u^\epsilon)'(C_T^\epsilon)^{-1}(Y_u^\epsilon)\right]^{1/2} \left[ (\dot Y_u^\epsilon)'(C_T^\epsilon)^{-1} (\dot Y_u^\epsilon)\right]^{1/2} du \nonumber \\
& \geq & -2 \left[ (\dot Y_T^\epsilon)'(C_T^\epsilon)^{-1} (\dot Y_T^\epsilon) \right]^{1/2} \left[ (Y_T^\epsilon)'(C_T^\epsilon)^{-1}(Y_T^\epsilon)\right]^{1/2} \nonumber \\ 
&& -2 \left[ (Y_T^\epsilon)'(C_T^\epsilon)^{-1} (Y_T^\epsilon) \right] \left[ \int_0^T (Y_u^\epsilon)'(C_T^\epsilon)^{-1}(Y_u^\epsilon) du + \int_0^T (\dot Y_u^\epsilon)'(C_T^\epsilon)^{-1} (\dot Y_u^\epsilon) du \right] \nonumber \\
&=& -2 \left[ (\dot Y_T^\epsilon)'(C_T^\epsilon)^{-1} (\dot Y_T^\epsilon) \right]^{1/2} \left[ (Y_T^\epsilon)'(C_T^\epsilon)^{-1}(Y_T^\epsilon)\right]^{1/2} \nonumber \\
&& -2 \left[ (Y_T^\epsilon)'(C_T^\epsilon)^{-1}(Y_T^\epsilon) \right] \left[ p + {\rm Tr} [(\dot C_T^\epsilon)(C_T^\epsilon)^{-1} ]\right]. \nn
\end{eqnarray}

Therefore,
\begin{eqnarray}
&& \frac{\partial}{\partial \epsilon} \ln
 \left[ (Y_T^\epsilon)'(C_T^\epsilon)^{-1} (Y_T^\epsilon) \right] \nn\\
&=& \left[ (Y_T^\epsilon)'(C_T^\epsilon)^{-1} (Y_T^\epsilon)\right]^{-1}
 \frac{\partial}{\partial \epsilon} \left[ (Y_T^\epsilon)'(C_T^\epsilon)^{-1}
 (Y_T^\epsilon) \right] \nn\\
& \geq & -2 \left[ \frac{(\dot Y_T^\epsilon)'(C_T^\epsilon)^{-1}
 (\dot Y_T^\epsilon)}{(Y_T^\epsilon)'(C_T^\epsilon)^{-1} (Y_T^\epsilon)}
 \right]^{1/2} -2 \left[ p + {\rm Tr}
 \left[ (\dot C_T^\epsilon)(C_T^\epsilon)^{-1} \right] \right] , \nn
\end{eqnarray}
which is bounded below (by a negative number possibly depending on $\epsilon$)
 uniformly for large values of $T$ by (\ref{eq3.4})
 and using the fact that both 
$(\dot Y_T^\epsilon)'(C_T^\epsilon)^{-1}(\dot Y_T^\epsilon)$
and
$(\dot Y_T^\epsilon)'(\dot C_T^\epsilon)^{-1}(\dot Y_T^\epsilon)$
have the same order and the latter has
the order as that of
 $(Y_T^\epsilon)'(C_T^\epsilon)^{-1}(Y_T^\epsilon)$.
Hence the proof of Lemma \ref{lm3.1}.  \qed 

{\bf Proof of Theorem \ref{th3.2}.}
Let $F^\epsilon = F - \epsilon I, \ \epsilon >0$. Since all eigenvalues of $F$
 are on the left half space, the real parts of all eigenvalues of
 $F^{\epsilon}$ are negative, i.e., $Y_t^\epsilon$ is a stable process.
 By Corollary \ref{cor3.1},
\begin{eqnarray}
{\lim}_{T \rightarrow \infty} (Y_T^\epsilon)'(C_T^\epsilon)^{-1}(Y_T^\epsilon) =0.  \nn
\end{eqnarray}
Let $f(\epsilon) = \ln (Y_T^\epsilon)(C_T^\epsilon)^{-1}(Y_T^\epsilon)$. 
Fix an $\epsilon_1 > 0$.
 $f$ is a continuous function on $[0,\epsilon_1]$ and is differentiable 
in $(0,\epsilon_1)$. Then
by the Mean Value Theorem, there exists an
 $\epsilon_0 \in (0,\epsilon_1 )$ such that 
\begin{eqnarray}
 f(\epsilon_1) -  f(0) = \epsilon_1 \frac{\partial}{\partial \epsilon} 
 f(\epsilon) \bigl|_{\epsilon = \epsilon_0}. \nn
\end{eqnarray}
That is,
\begin{eqnarray}
\label{eq3.5}
\frac{(Y_T^{\epsilon_1})'(C_T^{\epsilon_1})^{-1}(Y_T^{\epsilon_1})}{Y_T'C_T^{-1}Y_T} \geq 
\exp\left\{\epsilon_1 \frac{\partial}{\partial \epsilon} 
 f(\epsilon) \left|_{\epsilon = \epsilon_0} \right. \right\} ,
\end{eqnarray}
which is uniformly positive (i.e., bounded away from zero) for large values of
 $T$
 by Lemma \ref{lm3.1}.
Since
\begin{eqnarray}
{\lim}_{T \rightarrow \infty} (Y_T^\epsilon)'(C_T^\epsilon)^{-1}(Y_T^\epsilon) =0 \quad {\rm a.s.} \nn
\end{eqnarray}
 by (\ref{eq3.5})
\begin{eqnarray}
{\lim}_{T \rightarrow \infty} Y_T'C_T^{-1}Y_T=0 \quad {\rm a.s.} \nn
\end{eqnarray}
Hence the proof of Theorem \ref{th3.2}.  \qed 

\begin{corollary}
\label{cor3.2}
With the same assumptions and notations as in Lemma \ref{lm3.1},
\begin{eqnarray}
||C_T^{-1/2}|| = O(T^{-1/2}) \ \ \ {\rm a.s.} \nn
\end{eqnarray}
\end{corollary}

{\bf Proof.}
 Consider
\begin{eqnarray}
\frac{\partial}{\partial \epsilon} {\rm Tr} [(C_T^\epsilon)^{-1}]
 &=& -2 {\rm Tr} \left[ (C_T^\epsilon)^{-1} \int_0^T \left[ (Y_u^\epsilon )(\dot Y_u^\epsilon)' dt \right] (C_T^\epsilon )^{-1} \right] \nonumber \\
& \geq & -2{\rm Tr} (C_T^\epsilon)^{-1} \int_0^T \left[ (\dot Y_u^\epsilon)'(C_T^\epsilon)^{-1} (\dot Y_u^\epsilon)\right]^{1/2} \left[ (Y_u^\epsilon)'(C_T^\epsilon)^{-1} (Y_u^\epsilon) \right]^{1/2} du \nonumber \\
& \geq & - {\rm Tr} \left( C_T^{\epsilon} \right)^{-1} \left[ \int_0^T \left( \dot Y_u^\epsilon \right)'(C_T^\epsilon)^{-1}(\dot Y_u^\epsilon) du + \int_0^T (Y_u^\epsilon)'(C_T^\epsilon)^{-1}(Y_u^\epsilon)du \right] \nonumber \\
&=& -{\rm Tr}(C_T^\epsilon)^{-1} \left( {\rm Tr}\left[ (C_T^\epsilon)^{-1}\dot C_T^\epsilon \right] + {\rm Tr} \left[ (C_T^\epsilon)^{-1}(C_T^\epsilon) \right] \right) .  \nn
\end{eqnarray}
Hence 
$\frac{\partial}{\partial \epsilon} \ln {\rm Tr} [(C_T^\epsilon)^{-1}]
 \ge 
- \left( {\rm Tr}\left[ (C_T^\epsilon)^{-1}\dot C_T^\epsilon \right] +
 p \right) $
which is bounded below (by a negative number possibly depending on $\epsilon$)
uniformly for large values of $T$.
 Therefore, as in  (\ref{eq3.5}), by the Mean Value Theorem,
$
\frac{{\rm Tr}[(C_T^\epsilon)^{-1}]}{{\rm Tr}[(C_T)^{-1}]} 
$
is uniformly positive (i.e., bounded away from zero) for large values of $T$.
Since 
${\rm Tr}[(C_T^\epsilon)^{-1}] = O(T^{-1})$, we have
$ O({\rm Tr}[(C_T)^{-1}]) \le O({\rm Tr}[(C_T^\epsilon)^{-1}]) = O(T^{-1})$.
 Again, as
 for any positive definite matrix $K_T$, $O(||K_T||) = O({\rm Tr} (K_T))$,
 we obtain by Corollary \ref{cor3.1}(i),
 $||(C_T)^{-1/2}|| = ||(C_T^\epsilon)^{-1/2}|| =O(T^{-1/2})$. Hence the result.
\qed 

\begin{remark}
\label{rm3.2}
\end{remark}
It is clear from the arguments in the above corollary \ref{cor3.2} that,
for the eigenvalues of $F$ on the left half space,
$$
\frac{1}{T} \la_{\min}(C_T) = \frac{1}{T \la_{\max}(C_T^{-1})} > 0,
$$
almost surely, uniformly for large values of $T$, since
$T \la_{\max}(C_T^{-1})  = T ||C_T^{-1}|| \le T \ O(T^{-1}) = O(1) 
\quad {\rm a.s.}$

\section{General Ornstein-Uhlenbeck Processes}
\setcounter{equation}{0}
\setcounter{section}{5}
\label{sec.5}
\vspace{-0.4cm}
For the Ornstein-Uhlenbeck process defined in (\ref{eq1.1}) with
RANK condition (\ref{eq1.2}),
 we have considered the
  case in which all the eigenvalues of 
 $F$ have positive real parts
 and
 the  case in which all the eigenvalues of
$F$ have zero or negative real parts (i.e., zero eigenvalues, purely imaginary
and the eigenvalues with negative real parts).
 Now we combine these cases to discuss the {\it mixed model} in which
 $F$ can be decomposed into rational canonical form as follows:
\begin{eqnarray}
MF=GM =\left( \begin{array}{cccccc} G_0 & 0 \\ 0 & G_1 \end{array}\right)
\left( \begin{array}{cccccc} M_0 \\ M_1 \end{array}\right) =
\left( \begin{array}{cccccc} G_0M_0 \\ G_1 M_1 \end{array}\right) , \nn
\end{eqnarray}
where all the characteristic roots of $G_0$ lie in the right half space
and all the characteristic roots of $G_1$ lie
 on the left half space.
Let 
\begin{eqnarray}
\left( \begin{array}{cccccc} U_{0t} \\ U_{1t} \end{array}\right) = \left( \begin{array}{cccccc} M_0 \\ M_1 
 \end{array}\right)Y_t = MY_t . \nn
\end{eqnarray}
Then
\begin{eqnarray}
d\left( \begin{array}{cccccc} U_{0t} \\ U_{1t} \end{array}\right) = MdY_t  
=
 MFY_t dt + MAdW_t 
=
 \left( \begin{array}{cccccc} G_0 & 0 \\ 0 & G_1 \end{array}\right)
 \left( \begin{array}{cccccc} U_{0t} \\ U_{1t} \end{array}\right) dt + MAdW_t .
 \nn
\end{eqnarray}
Also,
\begin{eqnarray}
\label{eq5.1}
\left( \int_0^T dW_t Y_t' \right) M'&=& \left( \begin{array}{cccccc}
 \int_0^T dW_t U_{0t}' \\ \int_0^T dW_t U_{1t}'  \end{array}\right) \nn\\
{\rm and} \qquad 
M \left( \int_0^T Y_t Y_t' dt \right) M' &=& \left( \begin{array}{cccccc}
 \int_0^T U_{0t}
 U_{0t}' dt & \int_0^T U_{0t} U_{1t}' dt \\ \int_0^T U_{1t} U_{0t}' dt 
& \int_0^T U_{1t} U_{1t}' dt \end{array}\right) . \nn
\end{eqnarray}

Define, $C_{1T} =\int_0^T U_{1t}U_{1t}' dt$. 
We now derive the following result.

\begin{lemma}
\label{lm4.0}
 Suppose,
for the Ornstein-Uhlenbeck process defined in
(\ref{eq1.1}), the RANK condition (\ref{eq1.2}) holds. Then
\begin{eqnarray}
\label{eq4.8}
\Sigma_T^{-1} =\left[ D_T M \left( \int_0^T Y_t Y_t' dt \right) M' D_T' \right]^{-1} \rightarrow \left( \begin{array}{cccccc}
B^{-1} & 0 \\ 0 & I_{p_1} 
\end{array}\right) \quad {\rm a.s.}
\end{eqnarray}
where $B$ is defined in Section \ref{sec.2} (before (\ref{eq2.5})),
$I_{p_1}$ is a $p_1$-dimensional identity matrix
 and
\begin{eqnarray}
D_T = \left( \begin{array}{cccccc} e^{-G_0T} & 0 \\ 0 & C_{1T}^{-1/2}  \nn
\end{array}\right) .
\end{eqnarray}
\end{lemma}

{\bf Proof.}
Observing (\ref{eq5.1}), we obtain, by Theorem \ref{th2.1}, that
\begin{eqnarray}
{\lim}_{T \rightarrow \infty} e^{-G_0T}\left( \int_0^T U_{0t} U_{0t}' dt
 \right) e^{-G_0'T}
 =B \qquad {\rm  is \ positive \ definite \ a.s.} \nn
\end{eqnarray} 
Again,
 $(\Sigma_T)_{11} = C_{1T}^{-1/2} C_{1T} C_{1T}^{-1/2} = I_{p_1}$.
Hence, the proof is complete once we show
$e^{-G_0T} (\int_0^T U_{0t}U_{1t}' dt) C_{1T}^{-1/2} \to
 0_{p_0 \times p_1}$
 matrix almost surely, as $T \to \infty$.
Notice that, 
 by Corollary \ref{cor2.1},
\begin{eqnarray}
{\lim}_{T \rightarrow \infty} 
\int_0^T||e^{-G_0T} U_{0t} || dt < \infty {\rm  \qquad a.s.} \nn
\end{eqnarray}
and from 
Theorem \ref{th3.2}
\begin{eqnarray}
{\lim}_{T \rightarrow \infty} U_{1T}'C_{1T}^{-1}U_{1T}=0 \qquad {\rm a.s.} \nn
\end{eqnarray}
Therefore,
for all $\omega$ outside a null set, and for any given $\epsilon >0$, 
there exists $T_0(\omega ) > 0$
  such that for all
$t \geq T_0(\omega )$, $(U_{1t}'C_{1t}^{-1}U_{1t})^{1/2} < 
\epsilon/({\lim}_{T \rightarrow \infty} 
\int_0^T||e^{-G_0T} U_{0t}(\omega) || dt)$. 
Hence
\begin{eqnarray}
\left|\left| e^{-G_0T} (\int_0^T U_{0t}U_{1t}' dt) C_{1T}^{-1/2} \right|\right|
& \leq & \int_0^T ||e^{-G_0 T} U_{0t}U_{1t}'C_{1T}^{-1/2}||dt \nonumber \\
& \leq & \int_0^{T_0(\omega )} ||e^{-G_0 T}U_{0t}||
 \hspace{0.1cm} ||C_{1T}^{-1/2}U_{1t}||dt \nn\\
 &&  +
 \int_{T_0(\omega )}^T ||e^{-G_0T}U_{0t}||
 \hspace{0.1cm} ||C_{1T}^{-1/2}U_{1t}||dt  \nn
\end{eqnarray}
As $T \rightarrow \infty$, the first term goes to 0 since $T_0(\omega)$ is
 fixed.
The second term is less than $\epsilon$ by the choice of $T_0(\omega)$
 since
$C_{1t}$ is increasing in $t$ (in the sense that
$C_{1t_2} - C_{1t_1}$ is positive definite whenever $t_2 > t_1$)
 and
$||C_{1T}^{-1/2}U_{1t}||=(U_{1t}'C_{1T}^{-1}U_{1t})^{1/2}
 \leq (U_{1t}'C_{1t}^{-1}U_{1t})^{1/2}$.
As $\ep$ is arbitrary, the proof is complete.
\qed

We now observe that,
\begin{eqnarray}
\hat F_T -F &=& \left[ T^{-1/2} A \left( \int_0^T dW_t Y_t' \right) M'D_T' \right] \left[ D_T M \left( \int_0^T Y_t Y_t' dt \right) M'D_T' \right]^{-1} \nonumber \\
&& \times (T^{1/2}D_T M) \nn
\end{eqnarray}
and
\begin{eqnarray}
T^{-1/2}A \left( \int_0^T dW_t Y_t' \right) M'D_T' &=& \left( \begin{array}{cccccc}
T^{-1/2}e^{-G_0T}(\int_0^TU_{0t}dW_t')A' \\
T^{-1/2}C_{1T}^{-1/2}(\int_0^T U_{1t}dW_t')A' \\
\end{array}\right)' . \nn
\end{eqnarray}

The first term $T^{-1/2}A\left( \int_0^T dW_t U_{0t}' \right)
 e^{-G_0'T} =O(T^{-1/2})$ by Corollary \ref{cor2.1}(ii).
 To show the remaining terms converges to 0, we prove the following Theorem.
This theorem is in the spirit of Theorem 2.2 of Wei \cite{key15}, which
is presented for the discrete case.

\begin{theorem}
\label{th4.1}
\begin{eqnarray}
\frac{1}{\sqrt{T}}\left( \int_0^T dW_t U_{1t}' \right) C_{1T}^{-1/2} \rightarrow 0 \qquad {\rm  a.s. \ as \quad } T \to \infty . \nn
\end{eqnarray}
\end{theorem}

To prove Theorem \ref{th4.1}, we need the following lemmas.

\begin{lemma}
\label{lm4.1}
Fix $t_0 > 0$. Then,
\begin{eqnarray}
\int_{t_0}^T U_{1t}'C_{1t}^{-1}U_{1t} dt = O(\log T) \qquad {\rm  a.s. \ as \quad }
T \to \infty . \nn
\end{eqnarray}
\end{lemma}

{\bf Proof.} 
Notice that,
\begin{eqnarray}
\frac{d}{dt} \log |C_{1t} | &=& {\rm Tr} \left( C_{1t}^{-1} \frac{d}{dt}C_{1t} \right) \nonumber \\
&=& {\rm Tr} \left( C_{1t}^{-1} U_{1t}U_{1t}' \right) = 
U_{1t}'C_{1t}^{-1} U_{1t} , \nn
\end{eqnarray}
where $|C_{1t}|$ is the determinant of $C_{1t}$.
Observe that,
 $G_1$ can be further decomposed into a rational canonical form as follows:
\begin{eqnarray}
\left( \begin{array}{cccccc} M_{11} \\ M_{12} \\ M_{13} \end{array}\right) G_1
=
 \left( \begin{array}{cccccc} G_{11} & 0 & 0 \\ 0 & G_{12} & 0 \\
0 & 0 & G_{13} \end{array}\right) 
\left( \begin{array}{cccccc} M_{11} \\ M_{12} \\ M_{13} \end{array}\right) 
= \left( \begin{array}{cccccc} G_{11} M_{11} \\ G_{12} M_{12} \\ G_{13} M_{13} \end{array}\right) , \nn
\end{eqnarray}
where all the characteristic roots of $G_{11}$
 have negative real parts, those of
  $G_{12}$ are
 purely imaginary and those of $G_{13}$ are zero.
For $i, j = 1,  2, 3$, define
$C_{1tij} = \int_0^t U_{1is} U_{1js}' ds$, where
\begin{eqnarray}
\left( \begin{array}{cccccc} U_{11s} \\ U_{12s} \\ U_{13s} \end{array}\right) 
=
\left( \begin{array}{cccccc} M_{11} \\ M_{12} \\ M_{13} \end{array}\right) U_{1s} . \nn
\end{eqnarray}
Thus
$C_{1t} = \left(\left( C_{1tij} \right) \right)_{i,j=1,2,3}$,
and hence $|C_{1t}| \le |C_{1t11}| \ |C_{1t22}| \ |C_{1t33}|$.
Therefore, by Theorem \ref{th3.1} in Section \ref{sec.4}
 and Theorems \ref{th7.1} and \ref{th7.2} in the Appendix, one obtains
\begin{eqnarray}
\int_{t_0}^T U_{1t}'C_{1t}^{-1}U_{1t}dt = \log \frac{|C_{1T}|}{|C_{1t_0}|}=O(\log T) \qquad {\rm a.s. \ as \quad } T  \to \infty .  \nn
\end{eqnarray}
Hence, the proof.
\qed

\medskip
We observe that,
from Lemma \ref{lm4.1}, if we let
$g(T) = \int_{t_0}^T U_{1t}'C_{1t}^{-1}U_{1t} dt$, then $g(T) \uparrow \infty$
as $T \uparrow \infty$ almost surely.
Also, $E(\log |C_{1T}|) = E(\sum_{i} \log (\lambda_{i}(C_{1T})))$ $=
\sum_{i} E( \log (\lambda_{i}(C_{1T}))) \le$ $\sum_{i} \log(E(\lambda_{i}(C_{1T})))$ \\
 $\le p_1 \log (E(\lambda_{max}(C_{1T}))) \le p_1 \log \int_0^T E(\|U_{1t}\|^2)
 dt$. It is clear that, for the eigenvalues on the left half space,
$ E(\|U_{1t}\|^2)$ is at most $O(t^k)$, i.e.,
 it grows at most like a polynomial in $t$.
  Thus, $E(\log |C_{1T}|) = O(\log T)$ as well.
Hence, using integration by parts, we obtain,
\begin{eqnarray}
\label{eq4.19}
 E\left( \int_{t_0}^{\infty} \frac{U_{1t}'C_{1t}^{-1}U_{1t}}{t} dt \right)
< \infty .
\end{eqnarray}

\begin{lemma}
\label{lm4.2}
 Let $M_{1T} =\int_0^T dW_tU_{1t}'$. Then, under the hypothesis
of Theorem \ref{th4.1},
\begin{eqnarray}
\frac{1}{T^{1/2}} M_{1T} C_T^{-1/2}
 \rightarrow 0 \qquad {\rm  in \ probability.} \nn
\end{eqnarray}
\end{lemma}

{\bf Proof.}
Notice that $M_{1t}$ is a martingale with respect to the filtration 
$\{{\cal F}_t\}_{t\ge 0}$ where
${\cal F}_t = \sigma\{ W_s \ : \ 0 \le s \le t \}$.
Define 
 $N_{1T}=\int_{t_1}^T dW_t U_{1t}' =M_{1T}-M_{1t_1}$. Then,
 for $T>t_1$, $N_{1T}$ is also a martingale.
 Define $V_t =
 {\rm Tr}[C_{1t}^{-1}M_{1t}'M_{1t}]/ t$ and 
$\tilde V_t = {\rm Tr}[C_{1t}^{-1}N_{1t}'N_{1t}]/t$.
Since $\|\frac{1}{T^{1/2}} M_{1T} C_T^{-1/2}\|^2 \le V_T \le 2 \tilde V_T +
2 {\rm Tr}(C_{1T}^{-1}M_{1t_1}'M_{1t_1}]/ T$ and 
$ {\rm Tr}(C_{1T}^{-1}M_{1t_1}'M_{1t_1}]/ T \to 0$, almost surely, as
 $T \to \infty$, it is enough to show that $\tilde V_T \to 0$, in probability,
as $T \to \infty$ and this would be immediate once one shows $E(\tilde V_T)
 \to 0$ as $T \to \infty$.
Now use
 It\^{o}'s Lemma to get
\begin{eqnarray}
\label{eq4.1}
d\tilde V_t &=& \frac{\left[ {\rm Tr}\left( C_{1t}^{-1}d(N_{1t}'N_{1t})\right) + {\rm Tr} \left[ (\dot C_{1t}^{-1} ) N_{1t}'N_{1t} \right] dt \right]}{t}
 - \frac{\tilde V_t}{t} dt
\end{eqnarray}
where 
$\dot C_{1t}^{-1} = -C_{1t}^{-1} \left( \dot C_{1t} \right) C_{1t}^{-1} = -C_{1t}^{-1} U_{1t}U_{1t}'C_{1t}^{-1}$
which is non-positive definite. Thus, 
\\
${\rm Tr}\left[ \left( \dot C_{1t}^{-1} \right) N_{1t}'N_{1t} \right] =
-U_{1t}'C_{1t}^{-1}N_{1t}'N_{1t}C_{1t}^{-1}U_{1t} \leq 0$.
Therefore, by (\ref{eq4.1}) and applying the It\^{o}'s Lemma again,
 one obtains 
\begin{eqnarray}
\tilde V_T & \leq & \int_{t_1}^T {\rm Tr} \left( C_{1t}^{-1} d(N_{1t}'N_{1t}) \right)/t \nonumber \\
&=& \int_{t_1}^T {\rm Tr} \left( C_{1t}^{-1} \left[ (dN_{1t}')N_{1t} + N_{1t}'(dN_{1t}) + (dN_{1t}')(dN_{1t}) \right] \right)/t . \nn
\end{eqnarray}
Define $\tau_n = \inf \{ t>t_1 : |\tilde V_t | \geq n \}$, then
\begin{eqnarray}
E\tilde V_{T\wedge \tau_n} & \leq & E \int_{t_1}^{T \wedge \tau_n} {\rm Tr} \left( C_{1t}^{-1} (dN_{1t}')(dN_{1t}) \right)/t \\
&=& E \int_{t_1}^{T \wedge \tau_n} \frac{U_{1t}'C_{1t}^{-1} U_{1t} }{t}dt.  \nn
\end{eqnarray}
Since $V_{T \wedge \tau_n}$ and $U_{1t}'C_{1t}^{-1}U_{1t}$ are non-negative,
 by Fatou's Lemma and the Monotone Convergence Theorem,
\begin{eqnarray}
E \tilde V_T & \leq & E\int_{t_1}^T \frac{U_{1t}'C_{1t}^{-1}U_{1t}}{(\log t)^{1+\alpha}}dt.  \nn
\end{eqnarray}
Now, by the argument in (\ref{eq4.19}), one has
$\limsup_{\{T \to \infty\}}
E\tilde V_T 
\leq 
\alpha C t_1^{-\alpha} $.
As $t_1$ can be taken to be arbitrarily large,
we have the result.
\qed 

\begin{lemma}
\label{lm4.3}
Let $V_t = {\rm Tr} [C_{1t}^{-1}M_{1t}'M_{1t}]/t$. Then,
 with the same assumptions and notations as in Lemma \ref{lm4.2},
\begin{eqnarray}
\int_{t_1}^\infty E \left[ E(dV_t | {\cal F}_t ) \right]^+ < \infty . \nn
\end{eqnarray}
\end{lemma}

{\bf Proof.}
 Applying It\^{o}'s Lemma on $V_t$,
\begin{eqnarray}
dV_t &=& \frac{\left( {\rm Tr}\left[ C_{1t}^{-1} d(M_{1t}'M_{1t})\right] + {\rm Tr} \left[ \dot C_{1t}^{-1} (M_{1t}'M_{1t} )\right] dt \right)}{t}
 - \frac{V_t}{t} dt \nn
\end{eqnarray}
where
$
\dot C_{1t}^{-1} 
= -C_{1t}^{-1} \left(\dot C_{1t} \right) C_{1t}^{-1} = -C_{1t}^{-1} U_{1t}U_{1t}'C_{1t}^{-1}
$
and
$$
{\rm Tr} \left[ \dot C_{1t}^{-1}(M_{1t}'M_{1t}) \right]  
= -U_{1t}'C_{1t}^{-1}M_{1t}'M_{1t}C_{1t}^{-1}U_{1t} \leq 0. 
$$
Therefore,
\begin{eqnarray}
E(dV_t | {\cal F}_t) & \leq & E \left( \left[ {\rm Tr} (C_{1t}^{-1}d(M_{1t}'M_{1t}) ) \right] /t \ | \ {\cal F}_t \right) \nonumber \\
&=& E \left( \left[ {\rm Tr} \left(C_{1t}^{-1}[ (dM_{1t}')M_{1t} + M_{1t}'(dM_{1t}) + (dM_{1t}')(dM_{1t})] \right) \right]/t \ | \ {\cal F}_t \right) \nonumber \\
&=& E \left( \left[ {\rm Tr} \left( C_{1t}^{-1} (dM_{1t}')(dM_{1t}) \right) \right]/t \ | \ {\cal F}_t \right) \nonumber \\
&=& E \left( \frac{U_{1t}'C_{1t}^{-1}U_{1t}}{t} dt \ | \ {\cal F}_t \right) \nonumber \\
&=& \frac{U_{1t}'C_{1t}^{-1}U_{1t}}{t}dt. \nn
\end{eqnarray}
Thus,
\begin{eqnarray}
\left[E(dV_t | {\cal F}_t)\right]^+
 \leq \frac{U_{1t}'C_{1t}^{-1}U_{1t}}{t} dt. \nn
\end{eqnarray}
Since $U_{1t}'C_{1t}^{-1}U_{1t} \geq 0$, by Fubini's theorem and by
(\ref{eq4.19})
\begin{eqnarray}
\int_{t_1}^\infty E \left[E(dV_t | {\cal F}_t )\right]^+ 
=
E \int_{t_1}^\infty \left[E(dV_t | {\cal F}_t )\right]^+  
\le
E \int_{t_1}^\infty \frac{U_{1t}'C_{1t}^{-1}U_{1t}}{t} dt
<
 \infty. \nn
\end{eqnarray}
Hence, the proof. \qed 

{\bf Proof of Theorem \ref{th4.1}.}
 Define 
$A_\delta^{t_1,T} = \{ {\max}_{t_1 < t <T} V_t > \delta\}$ and 
$H_{t_1}=\{ V_{t_1} \leq \epsilon \}$ for any $\epsilon >0$. Then, using the
 Lenglart Inequality (cf. Karatzas and Shreve \cite{key3} p30 or Lenglart
 \cite{key16}),
\begin{eqnarray}
P\left(A_\delta^{t_1,T} \cap H_{t_1} \right)  \leq \frac{1}{\delta} EV_{t_1}I_{H_{t_1}} + \frac{1}{\delta} \int_{t_1}^T E \left( [E(dV_t | {\cal F}_t ) ]^+ I_{H_{t_1}} \right). \nn
\end{eqnarray}
Therefore,
\begin{eqnarray}
P\left(A_\delta ^{t_1,T} \right)
 & =&
 P \left( A_\delta ^{t_1 ,T} \cap H_{t_1}^c \right) 
 + P \left( A_\delta^{t_1,T} \cap H_{t_1} \right)
\nonumber \\
& \leq & 
P \left(H_{t_1}^c \right)
+ P \left( A_\delta ^{t_1,T} \cap H_{t_1} \right)
 \nonumber \\
& \leq & P \left( H_{t_1}^c \right) + \frac{1}{\delta } E V_{t_1}I_{H_{t_1}} + \frac{1}{\delta} \int_{t_1}^T E \left( \left[ E ( dV_t | {\cal F}_t ) \right]^+ I_{H_{t_1}} \right) \nonumber \\
& \leq & P \left( H_{t_1}^c \right) + \frac{\epsilon}{\delta} + \frac{1}{\delta} \int_{t_1}^\infty E \left[ E(dV_t | {\cal F}_t ) \right]^+ , \nn
\end{eqnarray}
which is finite since $\int_{t_1}^\infty E [E(dV_t | {\cal F}_t)]^+ < \infty$
 by Lemma \ref{lm4.3}.
Therefore, as $T \rightarrow \infty$,
\begin{eqnarray}
P \left( {\lim}_{ T \rightarrow \infty} A_\delta^{t_1,T} \right) &=& 
 {\lim}_{T \rightarrow \infty} P \left( A_\delta^{t_1,T} \right) \nonumber \\
& \leq & P\left( H_{t_1}^c \right) + \frac{\epsilon}{\delta} + \frac{1}{\delta} \int_{t_1}^\infty E \left[ E(dV_t | {\cal F}_t ) \right]^+. \nn
\end{eqnarray}
Thus,
\begin{eqnarray}
{\limsup}_{ t_1 \rightarrow \infty}\hspace{0.2cm} P \left( 
{\lim}_{T \rightarrow \infty} A_\delta^{t_1,T} \right) & \leq & \frac{\epsilon}{\delta} . \nn
\end{eqnarray}
Since this is true for all $\epsilon >0$,
\begin{eqnarray}
{\limsup}_{ t_1 \rightarrow \infty} \hspace{0.2cm} P \left(
 {\lim}_{ T \rightarrow \infty} A_\delta^{t_1,T} \right) = 0. \nn
\end{eqnarray}
This implies,
\begin{eqnarray}
\frac{1}{T^{1/2}}\left( \int_0^T dW_t U_{1t}' \right) C_{1T}^{-1/2} \rightarrow 0 \qquad {\rm  a.s.} \nn
\end{eqnarray}
Hence, the Theorem. 
\qed

{\bf Proof of Theorem \ref{th1.1}.}
\vspace{0.5cm}
From Lemma \ref{lm2.1}, we have
$||e^{-G_0T}|| =O(e^{-\la_0 T})$ and, from Corollary \ref{cor3.2}, we have
 $||C_{1T}^{-1/2}||=
O(T^{-1/2})$ almost surely, as $T \to \infty$.
 Thus,
\begin{eqnarray}
|| T^{1/2}D_T M|| = T^{1/2} ||M|| \left( ||e^{-G_0T}|| + ||C_{1T}^{-1/2}||
\right) = O(1) \quad {\rm a.s. \ as \ } T \to \infty . \nonumber 
\end{eqnarray}
Therefore, from 
(\ref{eq4.8}), Corollary \ref{cor2.1}(ii) and Theorem \ref{th4.1}, we have
$\lim_{T\to \infty} \hat{F}_T = F$ a.s.

\vspace{0.4cm}

To show that (\ref{eq1.13}) holds, we observe that, for
 the eigenvalues of $F$ in the
right half space
 (\ref{eq1.13})
 follows from Theorem \ref{th2.1} and, for the eigenvalues
of $F$ on the left half space 
 (\ref{eq1.13})
 follows from arguments in Corollary 
\ref{cor3.2} and Remark \ref{rm3.2}.
For the mixed model, we observe
\begin{eqnarray}
\left( \int_0^T Y_t Y_t' dt \right)^{-1} = D_T M \Sigma_T M'D_T' \nn
\end{eqnarray}
where $\lim_{T \rightarrow \infty} \Sigma_T$ is a.s. positive definite.
 Thus, by Lemma \ref{lm2.2},
\begin{eqnarray}
\lambda_{\max} \left[ \left( \int_0^T Y_t Y_t' dt \right)^{-1} \right] =
 O( \lambda_{\max} (D_T D_T')) =O(T^{-1}). \nn
\end{eqnarray}
Therefore, the Theorem follows.
\qed

\section{Asymptotic Efficiency}
\setcounter{equation}{0}
\setcounter{section}{6}
\label{sec.6}
\vspace{-0.4cm}

In this section we would like to show that our estimator for the
drift matrix $F$ is 
asymptotically efficient even if the underlying process is not
necessarily stationary (stable).
For matrix-valued estimator there several ways to define 
 asymptotic efficiency (see Barndorff-Nielson and S\/orensen \cite{barnsoren}, for details).

The result is already known in one-dimensional case and for vector-valued parameters
(e.g., \cite{key11, dtzkuto, kuto1, prakasa1} and references therein)
 when the
processes are not necessarily stationary. For multi-dimensional matrix-valued case, similar things can be
proved once the asymptotic efficiency is properly defined for the matrix
valued estimator.

Observe that, when $A A'$ is nonsingular, the log-likelihood of $F$,
(see \cite{key11}, pp. 213-214), on $[0, T]$ is defined by,
$L_{A}(F) = \int_{0}^{T} (Y_t' F' (AA')^{-1} dY_t) - 
(1/2)  \int_{0}^{T} (Y_t' F' (AA')^{-1} F Y_t) dt.$
Thus, 
$$
d L_{A}(F) = tr\left[ d F \left(\int_{0}^{T} Y_t dY_t'\right) (AA')^{-1} - 
 dF \left(\int_{0}^{T} Y_t Y_t' dt\right) F' (AA')^{-1} \right] .
$$
Therefore,
$d L_{A}(F)/dF = $ \\ 
$
\left(\int_{0}^{T} dY_t Y_t'\right)\left(\int_{0}^{T} Y_t Y_t' dt\right)^{-1}$.
When $A A'$ is not nonsingular, the log-likelihood of $F$ cannot
be written explicitly. Therefore, M.L.E. of $F$ could not be achieved. However,
we would show that the above estimator is asymptotically efficient under
the assumptions of the section 2.

We show that 
$E({\rm Tr}[(\hat{F}_{T} - F) E(C_{T}) (\hat{F}_{T} - F)'])^{1/2} = O(1)$
as $T \to \infty$.

Let $S_{T} = \left(\int_{0}^{T} A dW_t Y_t' \right)$,
and $C_{T} = \left(\int_{0}^{T} Y_t Y_t' dt \right)$ as before.  
We use ${\rm Tr}[(\hat{F}_{T} - F) E(C_{T}) (\hat{F}_{T} - F)'] =
{\rm Tr}[S_{T} C_{T}^{-1} E(C_{T}) C_{T}^{-1} S_{T}'] \le {\rm Tr}[S_{T} C_{T}^{-1} S_{T}'] {\rm Tr}[C_{T}^{-1} E(C_{T})]$ to prove the following result.

\bigskip
{\bf Proof of Theorem \ref{th1.2}}

\underline{Case 1:}
Eigenvalues of $F$ are in the positive half space.

\bigskip
Observe that,
${\rm Tr}(S_{T}C_{T}^{-1} S_{T}') = $ \\
 ${\rm Tr}(S_{T} e^{-F'T}(e^{-FT}C_{T}e^{-F'T})^{-1} e^{-FT} S_{T}')$.
Since $S_{T} e^{-F'T}$ is a Gaussian process and its mean zero
and variance $e^{-FT} E(C_{T}) e^{-F'T}$ converges 
(in fact, to $E(B)$) as $T \to \infty$,
 $S_{T} e^{-F'T}$ converges to a finite Gaussian random variable in
distribution. 
Also, from Theorem (\ref{th2.1}), as $T \to \infty$,
$e^{-FT}C_{T}e^{-F'T}$ converges almost surely to $B$
 (which is positive definite with probability one). 
Thus, we obtain
${\rm Tr}(S_{T} e^{-F'T}(e^{-FT}C_{T}e^{-F'T})^{-1} e^{-FT} S_{T}')$
converges in distribution to finite random variable with finite
expectation.

Now,
${\rm Tr}(C_{T}^{-1} E(C_{T})) = 
 {\rm Tr}((e^{-FT}C_{T}e^{-F'T})^{-1} (e^{-FT} E(C_{T}) e^{-F'T}))$,
and from Theorem (\ref{th2.1}), as $T \to \infty$,
$(e^{-FT} C_{T} e^{-F'T})^{-1}$ converges to $B^{-1}$ almost surely.
Also, 
$e^{-FT} E(C_{T}) e^{-F'T} = \int_{0}^{T} e^{-Ft} Y_0 Y_0' e^{-F't} dt
 + \int_{0}^{T} t e^{-Ft} A A' e^{-F't} dt $,
which is finite as $T \to \infty$.
Thus, it remains to show, as $T \to \infty$,
$E(e^{-FT} C_{T} e^{-F'T})^{-1}$ converges to $E(B^{-1})$ (which is finite).
First observe that,
$Z_t - Y_0  = \int_0^t e^{-Fs} A dW_s \ $ is a symmetric (Gaussian) martingale
and with $E|Z_t - Y_0|^2 \le E|Z - Y_0|^2 < \infty$. Thus
$M_{Z} = \max_{0 \le t < \infty} (Z_t - Y_0)$
exists and has finite expectation.
 Also, (by symmetry) 
$m_Z = \min_{0 \le t < \infty} (Z_t - Y_0)$ exists and has finite second
 moment.
For symmetric matrices $D_1$ and $D_2$, define,
 $D_1 \ge D_2$ if $D_1 - D_2$ is non-negative definite.
Therefore,
\begin{eqnarray}
e^{-FT} C_{T} e^{-F'T} 
&=&
\int_{0}^{T} e^{-Ft} Z_{T-t} Z_{T-t}' e^{-F't} dt \nn\\
&\ge &
 \int_{0}^{T} e^{-Ft} (m_Z + Y_0) (m_Z + Y_0)' e^{-F't} dt  \nn\\
&\ge &
 \int_{0}^{T_0} e^{-Ft} (m_Z + Y_0) (m_Z + Y_0)' e^{-F't} dt  \nn
\end{eqnarray}
for all $T \ge T_0$, for some $T_0 > 0$ ($T_0$ may be taken to be 1).
Thus,
$(e^{-FT} C_{T} e^{-F'T})^{-1} \le 
(\int_{0}^{T_0} e^{-Ft} (m_Z + Y_0) (m_Z + Y_0)' e^{-F't} dt)^{-1} 
$
for all $T \ge T_0$.
Since right hand side has finite expectation, using dominated convergence
type theorem deduce
$E(B^{-1}) = \lim_{T\to \infty}E(e^{-FT} C_{T} e^{-F'T})^{-1} 
\le E(\int_{0}^{T_0} e^{-Ft} (m_Z + Y_0) (m_Z + Y_0)' e^{-F't} dt)^{-1} .
$
Therefore,
$E({\rm Tr}(C_{T}^{-1} E(C_{T})))$ is finite
and hence
$E({\rm Tr}[(\hat{F}_{T} - F) E(C_{T}) (\hat{F}_{T} - F)'])^{1/2} = O(1)$.

\bigskip
\underline{Case 2:}
Eigenvalues of $F$ are on the left half space.

\bigskip
When all the eigenvalues have real parts negative,
by ergodic theorem,
$\lim_{T \to \infty} \ \frac{1}{T}C_{T} 
=  \int_0^{\infty} e^{Ft} A A' e^{F't} dt 
= \lim_{T \to \infty} \ E(\frac{1}{T}C_{T})$.
Thus, \\
$\lim_{T \to \infty} \ E({\rm Tr}(S_{T}C_{T}^{-1} S_{T}'))
= \lim_{T \to \infty} \ E({\rm Tr}(\frac{1}{T}S_{T}' S_{T}
 ( \int_0^{\infty} e^{Ft} A A' e^{F't} dt)^{-1}))
= p$, i.e., of O(1).
Also, 
$\lim_{T \to \infty} \ E({\rm Tr}(C_{T}^{-1} E(C_{T})))
= \lim_{T \to \infty} \ E({\rm Tr}((\frac{1}{T}C_{T})^{-1} E (\frac{1}{T}C_{T})))
= p$.
Therefore,
$E({\rm Tr}[(\hat{F}_{T} - F) E(C_{T}) (\hat{F}_{T} - F)'])^{1/2} = O(1)$.

\bigskip 
{\bf Zero and purely imaginary eigenvalues.}

When the eigenvalues are either all purely imaginary or all
zero,  
replace $F$ by $F - \ep I = F^{\ep}$,
 as it is done in Section \ref{sec.4},
get the result as above by ergodic theorem.

Now, as in Lemma \ref{lm3.1},
consider
\begin{eqnarray}
&& \frac{\partial}{\partial \epsilon} 
{\rm Tr}E((S_T^\epsilon)'(C_T^\epsilon)^{-1}(S_T^\epsilon)) \nonumber \\
&=& 2 {\rm Tr}E((\dot S_T^\epsilon)'(C_T^\epsilon)^{-1} S_T^\epsilon) + 
{\rm Tr}E((S_T^\epsilon)'\frac{\partial}{\partial \epsilon}(C_T^\epsilon)^{-1}
 S_T^\epsilon) \nonumber \\
&=& 2{\rm Tr}E((\dot S_T^\epsilon)'(C_T^\epsilon)^{-1}S_T^\epsilon) - 
{\rm Tr}E((S_T^\epsilon)'(C_T^\epsilon)^{-1} \left[ \frac{\partial}{\partial
 \epsilon} C_T^\epsilon \right] (C_T^\epsilon)^{-1} S_T^\epsilon)  \nonumber \\
& \geq & -2 E\left(\left[ {\rm Tr}((\dot S_T^\epsilon)'(C_T^\epsilon)^{-1}
 (\dot S_T^\epsilon ))\right]^{1/2} \left[{\rm Tr}((S_T^\epsilon)'
(C_T^\epsilon)^{-1} (S_T^\epsilon))\right]^{1/2}\right) \nonumber \\
&& - {\rm Tr}E\left[(S_T^\epsilon)'(C_T^\epsilon)^{-1}
 \left[ \int_0^T (Y_u^\epsilon)(\dot Y_u^\epsilon)'du + 
\int_0^T (\dot Y_u^\epsilon)(Y_u^\epsilon)'du \right] 
(C_T^\epsilon)^{-1}S_T^\epsilon\right] \nonumber \\
& \geq & -2 \left(E\left[{\rm Tr}((\dot S_T^\epsilon)'(C_T^\epsilon)^{-1}
 (\dot S_T^\epsilon)) \right]\right)^{1/2} \left(E\left[ {\rm Tr}((S_T^\epsilon)'
(C_T^\epsilon)^{-1}(S_T^\epsilon)) \right]\right)^{1/2}  \nonumber \\
&& -2 E\left(\left[ {\rm Tr}((S_T^\epsilon)'(C_T^\epsilon)^{-1}(S_T^\epsilon))\right]
 \int_0^T \left[ (Y_u^\epsilon)'(C_T^\epsilon)^{-1}(Y_u^\epsilon)\right]^{1/2}
 \left[ (\dot Y_u^\epsilon)'(C_T^\epsilon)^{-1} 
(\dot Y_u^\epsilon)\right]^{1/2} du \right) \nonumber \\
& \geq & -2 \left(E\left[ {\rm Tr}((\dot S_T^\epsilon)'(C_T^\epsilon)^{-1}
 (\dot S_T^\epsilon)) \right]\right)^{1/2} \left(E\left[ {\rm Tr}((S_T^\epsilon)'
(C_T^\epsilon)^{-1}(S_T^\epsilon)) \right]\right)^{1/2} \nonumber \\ 
&& - E\left(\left[ {\rm Tr}((S_T^\epsilon)'(C_T^\epsilon)^{-1} (S_T^\epsilon)) \right]
 \left[ \int_0^T (Y_u^\epsilon)'(C_T^\epsilon)^{-1}(Y_u^\epsilon) du + 
\int_0^T (\dot Y_u^\epsilon)'(C_T^\epsilon)^{-1}
 (\dot Y_u^\epsilon) du \right] \right) \nonumber \\
&=& -2 \left(E\left[ {\rm Tr}((\dot S_T^\epsilon)'(C_T^\epsilon)^{-1}
 (\dot S_T^\epsilon)) \right]\right)^{1/2} \left(E\left[ {\rm Tr}((S_T^\epsilon)'
(C_T^\epsilon)^{-1}(S_T^\epsilon)) \right]\right)^{1/2} \nonumber \\
&& - E\left(\left[ {\rm Tr}((S_T^\epsilon)'(C_T^\epsilon)^{-1}(S_T^\epsilon)) \right] 
\left[ p + {\rm Tr} [(\dot C_T^\epsilon)(C_T^\epsilon)^{-1} ]\right]\right). \nn
\end{eqnarray}

Therefore,
\begin{eqnarray}
&& \frac{\partial}{\partial \epsilon} \ln
 E\left[ {\rm Tr}((S_T^\epsilon)'(C_T^\epsilon)^{-1} (S_T^\epsilon)) \right] \nn\\
&=& \left[E {\rm Tr}((S_T^\epsilon)'(C_T^\epsilon)^{-1} (S_T^\epsilon)) \right]^{-1}
 \frac{\partial}{\partial \epsilon} 
E\left[ {\rm Tr}((S_T^\epsilon)'(C_T^\epsilon)^{-1} (S_T^\epsilon)) \right] \nn\\
& \geq & -2 \left[ \frac{E{\rm Tr}((\dot S_T^\epsilon)'(C_T^\epsilon)^{-1}
 (\dot S_T^\epsilon))}{E{\rm Tr}((S_T^\epsilon)'(C_T^\epsilon)^{-1} (S_T^\epsilon))}
 \right]^{1/2} -  \frac{E\left(\left[ {\rm Tr}((S_T^\epsilon)'(C_T^\epsilon)^{-1} (S_T^\epsilon)) \right]\left[ p + {\rm Tr}
 \left[ (\dot C_T^\epsilon)(C_T^\epsilon)^{-1} \right] \right]\right)}
{E\left[ {\rm Tr}((S_T^\epsilon)'(C_T^\epsilon)^{-1} (S_T^\epsilon)) \right]} , \nn
\end{eqnarray}
which is bounded below (by a negative number possibly depending on $\epsilon$)
 uniformly for large values of $T$ by (\ref{eq3.4})
 and using the fact that both 
${\rm Tr}E((\dot S_T^\epsilon)'(C_T^\epsilon)^{-1}(\dot S_T^\epsilon))$
and
${\rm Tr}E((\dot S_T^\epsilon)'(\dot C_T^\epsilon)^{-1}(\dot S_T^\epsilon))$
have the same order and the latter has
the order as that of
 ${\rm Tr}E((S_T^\epsilon)'(C_T^\epsilon)^{-1}(S_T^\epsilon))$. 

Now as in the argument 
in consistency part,
since all eigenvalues of $F$
 are on the left half space, the real parts of all eigenvalues of
 $F^{\epsilon}$ are negative, i.e., $Y_t^\epsilon$ is a stable process
and
\begin{eqnarray}
{\lim}_{T \rightarrow \infty} {\rm Tr}E((S_T^\epsilon)'(C_T^\epsilon)^{-1}(S_T^\epsilon)) = O(1).  \nn
\end{eqnarray}

Similarly, to get a upper bound, consider
\begin{eqnarray}
&& \frac{\partial}{\partial \epsilon} 
{\rm Tr}E((S_T^\epsilon)'(C_T^\epsilon)^{-1}(S_T^\epsilon)) \nonumber \\
&=& 2{\rm Tr}E((\dot S_T^\epsilon)'(C_T^\epsilon)^{-1}S_T^\epsilon) - 
{\rm Tr}E((S_T^\epsilon)'(C_T^\epsilon)^{-1} \left[ \frac{\partial}{\partial
 \epsilon} C_T^\epsilon \right] (C_T^\epsilon)^{-1} S_T^\epsilon)  \nonumber \\
& \leq & 2 E\left(\left[ {\rm Tr}((\dot S_T^\epsilon)'(C_T^\epsilon)^{-1}
 (\dot S_T^\epsilon ))\right]^{1/2} \left[{\rm Tr}((S_T^\epsilon)'
(C_T^\epsilon)^{-1} (S_T^\epsilon))\right]^{1/2}\right) \nonumber \\
&& + {\rm Tr}E\left[(S_T^\epsilon)'(C_T^\epsilon)^{-1}
 \left[ \int_0^T (Y_u^\epsilon)(\dot Y_u^\epsilon)'du + 
\int_0^T (\dot Y_u^\epsilon)(Y_u^\epsilon)'du \right] 
(C_T^\epsilon)^{-1}S_T^\epsilon\right] \nonumber \\
& \leq & 2 \left(E\left[{\rm Tr}((\dot S_T^\epsilon)'(C_T^\epsilon)^{-1}
 (\dot S_T^\epsilon)) \right]\right)^{1/2} \left(E\left[ {\rm Tr}((S_T^\epsilon)'
(C_T^\epsilon)^{-1}(S_T^\epsilon)) \right]\right)^{1/2}  \nonumber \\
&& + 2 E\left(\left[ {\rm Tr}((S_T^\epsilon)'(C_T^\epsilon)^{-1}(S_T^\epsilon))\right]
 \int_0^T \left[ (Y_u^\epsilon)'(C_T^\epsilon)^{-1}(Y_u^\epsilon)\right]^{1/2}
 \left[ (\dot Y_u^\epsilon)'(C_T^\epsilon)^{-1} 
(\dot Y_u^\epsilon)\right]^{1/2} du \right) \nonumber \\
& \leq & 2 \left(E\left[ {\rm Tr}((\dot S_T^\epsilon)'(C_T^\epsilon)^{-1}
 (\dot S_T^\epsilon)) \right]\right)^{1/2} \left(E\left[ {\rm Tr}((S_T^\epsilon)'
(C_T^\epsilon)^{-1}(S_T^\epsilon)) \right]\right)^{1/2} \nonumber \\ 
&& + E\left(\left[ {\rm Tr}((S_T^\epsilon)'(C_T^\epsilon)^{-1} (S_T^\epsilon)) \right]
 \left[ \int_0^T (Y_u^\epsilon)'(C_T^\epsilon)^{-1}(Y_u^\epsilon) du + 
\int_0^T (\dot Y_u^\epsilon)'(C_T^\epsilon)^{-1}
 (\dot Y_u^\epsilon) du \right] \right) \nonumber \\
&=& 2 \left(E\left[ {\rm Tr}((\dot S_T^\epsilon)'(C_T^\epsilon)^{-1}
 (\dot S_T^\epsilon)) \right]\right)^{1/2} \left(E\left[ {\rm Tr}((S_T^\epsilon)'
(C_T^\epsilon)^{-1}(S_T^\epsilon)) \right]\right)^{1/2} \nonumber \\
&& + E\left(\left[ {\rm Tr}((S_T^\epsilon)'(C_T^\epsilon)^{-1}(S_T^\epsilon)) \right] 
\left[ p + {\rm Tr} [(\dot C_T^\epsilon)(C_T^\epsilon)^{-1} ]\right]\right). \nn
\end{eqnarray}

Therefore,
\begin{eqnarray}
&& \frac{\partial}{\partial \epsilon} \ln
 E\left[ {\rm Tr}((S_T^\epsilon)'(C_T^\epsilon)^{-1} (S_T^\epsilon)) \right] \nn\\
&=& \left[E {\rm Tr}((S_T^\epsilon)'(C_T^\epsilon)^{-1} (S_T^\epsilon)) \right]^{-1}
 \frac{\partial}{\partial \epsilon} 
E\left[ {\rm Tr}((S_T^\epsilon)'(C_T^\epsilon)^{-1} (S_T^\epsilon)) \right] \nn\\
& \leq & 2 \left[ \frac{E{\rm Tr}((\dot S_T^\epsilon)'(C_T^\epsilon)^{-1}
 (\dot S_T^\epsilon))}{E{\rm Tr}((S_T^\epsilon)'(C_T^\epsilon)^{-1} (S_T^\epsilon))}
 \right]^{1/2} 
 +  \frac{E\left(\left[ {\rm Tr}((S_T^\epsilon)'(C_T^\epsilon)^{-1} (S_T^\epsilon)) \right]\left[ p + {\rm Tr}
 \left[ (\dot C_T^\epsilon)(C_T^\epsilon)^{-1} \right] \right]\right)}
{E\left[ {\rm Tr}((S_T^\epsilon)'(C_T^\epsilon)^{-1} (S_T^\epsilon)) \right]} , \nn
\end{eqnarray}
which is bounded above (by a positive number possibly depending on $\epsilon$)
 uniformly for large values of $T$ by (\ref{eq3.4}).

As in the proof of Theorem \ref{th3.2},
let $f(\epsilon) = \ln {\rm Tr} E((S_T^\epsilon)(C_T^\epsilon)^{-1}(S_T^\epsilon))$. 
Fix an $\epsilon_1 > 0$.
 $f$ is a continuous function on $[0,\epsilon_1]$ and is differentiable 
in $(0,\epsilon_1)$. Then
by the Mean Value Theorem, there exists an
 $\epsilon_0 \in (0,\epsilon_1 )$ such that 
\begin{eqnarray}
 f(\epsilon_1) -  f(0) = \epsilon_1 \frac{\partial}{\partial \epsilon} 
 f(\epsilon) \bigl|_{\epsilon = \epsilon_0}. \nn
\end{eqnarray}
That is,
\begin{eqnarray}
\label{eq6.5}
\frac{{\rm Tr}E((S_T^{\epsilon_1})'(C_T^{\epsilon_1})^{-1}(S_T^{\epsilon_1}))}{{\rm Tr}E(S_T'C_T^{-1}S_T)} = 
\exp\left\{\epsilon_1 \frac{\partial}{\partial \epsilon} 
 f(\epsilon) \left|_{\epsilon = \epsilon_0} \right. \right\} ,
\end{eqnarray}
which is uniformly bounded and positive (i.e., bounded away from zero and infinity) for large values of
 $T$
as argued above.
Since
\begin{eqnarray}
{\lim}_{T \rightarrow \infty} {\rm Tr} E((S_T^\epsilon)'(C_T^\epsilon)^{-1}(S_T^\epsilon)) = O(1) 
. \nn
\end{eqnarray}
 by (\ref{eq6.5})
\begin{eqnarray}
{\lim}_{T \rightarrow \infty} {\rm Tr} E(S_T'C_T^{-1}S_T) = O(1)  
. \nn
\end{eqnarray} 

Mimicking the above argument, find
\begin{eqnarray}
&& \frac{\partial}{\partial \epsilon} 
{\rm Tr}\left(E((C_T^\epsilon)^{-1})E(C_T^\epsilon)\right) \nonumber \\
&=& - {\rm Tr}\left(E\left((C_T^\epsilon)^{-1} 
 \left[ \int_0^T (Y_u^\epsilon)(\dot Y_u^\epsilon)'du + 
\int_0^T (\dot Y_u^\epsilon)(Y_u^\epsilon)'du \right] 
 (C_T^\epsilon)^{-1}\right)E(C_T^\epsilon)\right) \nn\\
&& + {\rm Tr}\left(E((C_T^\epsilon)^{-1})
E \left[ \int_0^T (Y_u^\epsilon)(\dot Y_u^\epsilon)'du + 
\int_0^T (\dot Y_u^\epsilon)(Y_u^\epsilon)'du \right] \right) \nonumber \\
& \geq &
 -2 E\left(\left[ {\rm Tr}((C_T^\epsilon)^{-1}E(C_T^\epsilon))\right]
 \int_0^T \left[ (Y_u^\epsilon)'(C_T^\epsilon)^{-1}(Y_u^\epsilon)\right]^{1/2}
 \left[ (\dot Y_u^\epsilon)'(C_T^\epsilon)^{-1} 
(\dot Y_u^\epsilon)\right]^{1/2} du \right)  \nonumber \\
&&
 -2 E\left(
 \int_0^T \left[ (Y_u^\epsilon)'(E((C_T^\epsilon)^{-1}))
(Y_u^\epsilon)\right]^{1/2}
 \left[ (\dot Y_u^\epsilon)'(E((C_T^\epsilon)^{-1}))
(\dot Y_u^\epsilon)\right]^{1/2} du \right) \nonumber \\
& \geq & 
- E\left(\left[ {\rm Tr}((C_T^\epsilon)^{-1}E(C_T^\epsilon))\right]
\left[ p + {\rm Tr} [(\dot C_T^\epsilon)(C_T^\epsilon)^{-1} ]\right]\right)
 \nn\\
&& 
- 2 E\left(\left[ {\rm Tr}(E((C_T^\epsilon)^{-1}))(C_T^\epsilon)\right]^{1/2}
\left[ {\rm Tr}(E((C_T^\epsilon)^{-1}))(\dot C_T^\epsilon)\right]^{1/2}\right)
 \nn\\
& \geq & 
- E\left(\left[ {\rm Tr}((C_T^\epsilon)^{-1}E(C_T^\epsilon))\right]
\left[ p + {\rm Tr} [(\dot C_T^\epsilon)(C_T^\epsilon)^{-1} ]\right]\right)
 \nn\\
&& 
- 2
 \left({\rm Tr}\left[E((C_T^\epsilon)^{-1})E(C_T^\epsilon)\right]\right)^{1/2}
 \left({\rm Tr}\left[E((C_T^\epsilon)^{-1})E(\dot C_T^\epsilon)\right]\right)^{1/2}
 \nn
\end{eqnarray}

Therefore,
\begin{eqnarray}
&& \frac{\partial}{\partial \epsilon} \ln
 {\rm Tr}(E((C_T^\epsilon)^{-1}) E(C_T^\epsilon)) \nn\\
&=&
 \left[{\rm Tr}(E((C_T^\epsilon)^{-1}) E(C_T^\epsilon))\right]^{-1} 
 \frac{\partial}{\partial \epsilon}
 {\rm Tr}(E((C_T^\epsilon)^{-1}) E(C_T^\epsilon)) \nn\\
& \geq &
- \frac{E\left(\left[ {\rm Tr}((C_T^\epsilon)^{-1}E(C_T^\epsilon))\right]
\left[ p + {\rm Tr} [(\dot C_T^\epsilon)(C_T^\epsilon)^{-1} ]\right]\right)}
{{\rm Tr}(E((C_T^\epsilon)^{-1}) E(C_T^\epsilon))}
 -2 \left[\frac{{\rm Tr}(E((C_T^\epsilon)^{-1})E(\dot C_T^\epsilon))}
{{\rm Tr}(E((C_T^\epsilon)^{-1}) E(C_T^\epsilon))}\right]^{1/2} 
 , \nn
\end{eqnarray}
which is bounded below (by a negative number possibly depending on $\epsilon$)
 uniformly for large values of $T$ by (\ref{eq3.4})
 and using the fact that both 
${\rm Tr}(E((C_T^\epsilon)^{-1})E(\dot C_T^\epsilon))$
and
${\rm Tr}(E((C_T^\epsilon)^{-1})E(C_T^\epsilon))$
have the same order.

Similary, to get an upper bound, consider
\begin{eqnarray}
&& \frac{\partial}{\partial \epsilon} 
{\rm Tr}\left(E((C_T^\epsilon)^{-1})E(C_T^\epsilon)\right) \nonumber \\
&=& - {\rm Tr}\left(E\left((C_T^\epsilon)^{-1} 
 \left[ \int_0^T (Y_u^\epsilon)(\dot Y_u^\epsilon)'du + 
\int_0^T (\dot Y_u^\epsilon)(Y_u^\epsilon)'du \right] 
 (C_T^\epsilon)^{-1}\right)E(C_T^\epsilon)\right) \nn\\
&& + {\rm Tr}\left(E((C_T^\epsilon)^{-1})
E \left[ \int_0^T (Y_u^\epsilon)(\dot Y_u^\epsilon)'du + 
\int_0^T (\dot Y_u^\epsilon)(Y_u^\epsilon)'du \right] \right) \nonumber \\
& \leq &
  2 E\left(\left[ {\rm Tr}((C_T^\epsilon)^{-1}E(C_T^\epsilon))\right]
 \int_0^T \left[ (Y_u^\epsilon)'(C_T^\epsilon)^{-1}(Y_u^\epsilon)\right]^{1/2}
 \left[ (\dot Y_u^\epsilon)'(C_T^\epsilon)^{-1} 
(\dot Y_u^\epsilon)\right]^{1/2} du \right)  \nonumber \\
&&
 + 2 E\left(
 \int_0^T \left[ (Y_u^\epsilon)'(E((C_T^\epsilon)^{-1}))
(Y_u^\epsilon)\right]^{1/2}
 \left[ (\dot Y_u^\epsilon)'(E((C_T^\epsilon)^{-1}))
(\dot Y_u^\epsilon)\right]^{1/2} du \right) \nonumber \\
& \leq & 
 E\left(\left[ {\rm Tr}((C_T^\epsilon)^{-1}E(C_T^\epsilon))\right]
\left[ p + {\rm Tr} [(\dot C_T^\epsilon)(C_T^\epsilon)^{-1} ]\right]\right)
 \nn\\
&& 
+ 2 E\left(\left[ {\rm Tr}(E((C_T^\epsilon)^{-1}))(C_T^\epsilon)\right]^{1/2}
\left[ {\rm Tr}(E((C_T^\epsilon)^{-1}))(\dot C_T^\epsilon)\right]^{1/2}\right)
 \nn\\
& \leq & 
 E\left(\left[ {\rm Tr}((C_T^\epsilon)^{-1}E(C_T^\epsilon))\right]
\left[ p + {\rm Tr} [(\dot C_T^\epsilon)(C_T^\epsilon)^{-1} ]\right]\right)
 \nn\\
&& 
+ 2
 \left({\rm Tr}\left[E((C_T^\epsilon)^{-1})E(C_T^\epsilon)\right]\right)^{1/2}
 \left({\rm Tr}\left[E((C_T^\epsilon)^{-1})E(\dot C_T^\epsilon)\right]\right)^{1/2}
 \nn
\end{eqnarray}

Therefore,
\begin{eqnarray}
&& \frac{\partial}{\partial \epsilon} \ln
 {\rm Tr}(E((C_T^\epsilon)^{-1}) E(C_T^\epsilon)) \nn\\
&=&
 \left[{\rm Tr}(E((C_T^\epsilon)^{-1}) E(C_T^\epsilon))\right]^{-1} 
 \frac{\partial}{\partial \epsilon}
 {\rm Tr}(E((C_T^\epsilon)^{-1}) E(C_T^\epsilon)) \nn\\
& \leq &
 \frac{E\left(\left[ {\rm Tr}((C_T^\epsilon)^{-1}E(C_T^\epsilon))\right]
\left[ p + {\rm Tr} [(\dot C_T^\epsilon)(C_T^\epsilon)^{-1} ]\right]\right)}
{{\rm Tr}(E((C_T^\epsilon)^{-1}) E(C_T^\epsilon))}
 + 2 \left[\frac{{\rm Tr}(E((C_T^\epsilon)^{-1})E(\dot C_T^\epsilon))}
{{\rm Tr}(E((C_T^\epsilon)^{-1}) E(C_T^\epsilon))}\right]^{1/2} 
 , \nn
\end{eqnarray}
which is bounded above (by a positive number possibly depending on $\epsilon$)
 uniformly for large values of $T$ by (\ref{eq3.4}).
 
Thus, using the similar argument as in (\ref{eq6.5})
we show,
since 
${\lim}_{T \rightarrow \infty} {\rm Tr}(E((C_T^{\epsilon_1})^{-1})E(C_T^{\epsilon_1})) = O(1)
$,
$ \ \
{\lim}_{T \rightarrow \infty} {\rm Tr}(E(C_T^{-1})E(C_T)) = O(1) .
$
Hence, for eigenvalues of $F$ on the left half space, we prove that
$E({\rm Tr}[(\hat{F}_{T} - F) E(C_{T}) (\hat{F}_{T} - F)'])^{1/2} = O(1)$.

\bigskip
\underline{Case 3:} Mixed model.

In this case, use the decomposition of $F$ as in Section \ref{sec.5},
to decompose $Y_t'M' = (U_{0t}', U_{1t}')$.
Then,  one gets,
\begin{eqnarray}
tr(S_T {C_T}^{-1} S_T') 
& = &
tr(S_T M'D_T' (D_T M C_T M' D_T')^{-1} D_T M S_T') \nn\\
& \le & tr(S_T M' D_T' D_T M S_T') tr(D_T M C_T M' D_T')^{-1} \nn\\
& \le & (tr(S_{0T} e^{-G_0' T} e^{-G_0 T} S_{0T}') + tr(S_{1T} {C_{1T}}^{-1}
 S_{1T}')) tr(D_T M C_T M' D_T')^{-1} . \nn
\end{eqnarray}
Since for a symmetric invertible partition matrix,
$$
K = \left[\begin{array}{ll}
 E & F \\
 F' & H \\
\end{array}
\right]
$$
with $E$ and $H$ invertible, $tr(K) = tr(E - FH^{-1}F')^{-1} +
 tr(H - F'E^{-1}F)^{-1}$.
Taking $E = e^{- G_0 T} C_{0T} e^{- G_0' T}$,
$F = e^{- G_0 T} \int_0^T U_{0t} U_{1t}' dt C_{1T}^{-1/2} $
and $H = I$, i.e., identity matrix of order $p_1$.
Since $F$ converging to zero almost surely by the proof of Lemma \ref{lm4.0}
and by the same lemma $E$ converges to $B$ almost surely, one obtains
$tr(D_T M C_T M' D_T')^{-1} \to tr(B^{-1}) + p_1$ almost surely, as $T \to 
\infty$.
Therefore,
\bea
&& 
E\left[ (tr( e^{-G_0 T} S_{0T}' S_{0T} e^{-G_0' T}) + tr(S_{1T} {C_{1T}}^{-1}
 S_{1T}')) tr(D_T M C_T M' D_T')^{-1}\right]^{1/2} \nn\\
&\le & E(tr( e^{-G_0 T} S_{0T}' S_{0T} e^{-G_0' T}))
E(tr(D_T M C_T M' D_T')^{-1})  \nn\\
&& + 
E(tr(S_{1T} {C_{1T}}^{-1} S_{1T}'))
 E(tr(D_T M C_T M' D_T')^{-1}) \nn\\
&=& O(1)
\eea
by the case 1, and case 2. 
Similarly, \\
$tr((D_T M C_T M' D_T')^{-1} D_T E(M C_T M') D_T' ) 
\le tr((D_T M C_T M' D_T')^{-1}) tr(D_T E(M C_T M') D_T')
$
and 
$tr(D_T E(M C_T M') D_T') = tr(e^{- G_0 T} E(C_{0T}) e^{- G_0' T}) 
+ tr(C_{1T}^{-1} E(C_{1T}))$
expectation of which is finite by case 1 and case 2.
Therefore one proves, for the mixed model,
$E({\rm Tr}[(\hat{F}_{T} - F) E(C_{T}) (\hat{F}_{T} - F)'])^{1/2} = O(1)$.
\qed

\bigskip
{\bf Concluding remarks and discussion}

It is easy to see that the state space equation of
the general continuous autoregressive process (CAR(p)) of the form
$dX_{t}^{p-1} =  \alpha_p X_t + \alpha_{p-1} X_t^{1} + \cdots + 
\alpha_1 X_t^{p-1} + \sigma d W_t$
is a special case of multidimensional 
OU processes
where 
\begin{eqnarray}
F = \left(\begin{array}{cc}
{\bf 0}_{(p-1)\times 1} & {\bf I}_{p-1}\\
\alpha_p & \cdots \ \ \alpha_1
\end{array} \right), \hspace{0.25in}
A=\left(0,\cdots,0, \sigma \right)'  \nn
\end{eqnarray}
with $\alpha_i$ real numbers, $\sigma > 0$ and $W_t$ a one-dimensional
Browian motion. Clearly, $A$ is not singular. However,
the {\em RANK} condition (a) holds for this $F$ and $A$
 and, the condition (b') holds for this $F$.
Hence, from our result, the consistency and the asymptotic efficiency
of the $\hat{F}$
 of general
CAR(p)
follows.

It is important to observe that
 this estimation procedure may be the first step
 in developing a test of zero roots of some $F$,
which is necessary to determine 
whether univariate processes are co-integrated.
Also, if one needs to develop a test to determine whether the model for $Y_t$ 
is stationary, it is often enough to test whether all eigenvalues of $F$
have negative real parts against the alternative that some of them have
zero real parts.
 Therefore, one need not often worry about the assumption (b) or (b')
for testing stationarity.
Thus,
 a related
 question arises on, whether any 
  Asymptotically Mixed Normality 
property holds for the estimator $\hat F_T$, i.e., whether
$(\int_0^T Y_t Y_t' dt)^{1/2} (\hat F_T - F)$
 follows asymptotically Normal, so that we could compute approximate confidence
interval for the above testing procedures for the necessary parameters in $F$.
 As far as we know, these results are still unknown. Investigating the 
Asymptotically Mixed Normality 
property may be an important future direction to consider. One can look into
LAMN property as well.

Besides, when the drift coefficient matrix depends
on an unknown discrete paratmeter $\theta$ which 
follows a Markov chain (that helps the process to switch regimes),
 finding a consistent and asymptotically efficient estimator becomes important. 
Above questions can be asked in that
setup as well.

In applications, we almost always use discrete sampled data.
Similar questions can be asked for this model, when the data sampled are in deterministic (equal or unequal) time interval or in random interval. 
That can also be a focus of the future direction.

\section{Appendix}
\setcounter{equation}{0}
\setcounter{section}{7}
\label{sec.7}

\subsection{ Purely Imaginary Eigenvalues}
\vspace{-0.4cm}
 In this Section, we study the asymptotic behavior of
OU
 processes when the drift matrix $F$ only contains purely
 imaginary eigenvalues. The main results are summarized in the following:

\begin{theorem}
\label{th7.1}
 Suppose for the Ornstein-Uhlenbeck process defined in (\ref{eq1.1}),
the RANK condition (\ref{eq1.2}) holds and all the eigenvalues of $F$
 are purely imaginary. Let $2\rho$ be the dimension of the largest block of
 the rational canonical form of $F$ as 
defined in
 Section \ref{sec.2}
(see the Example). Then
\begin{eqnarray}
||Y_T|| = \left\{ \begin{array}{cccccc} O(T^{1/2}  \sqrt{ \ln \ln T}) & {\rm  a.s. \ \ \ if \ \  }
\rho =1  \\
O(T^{2\rho -5/2} \sqrt{\ln \ln T}) & {\rm  a.s. \ \ \  if \ \  }\rho \geq 2.
\end{array} \right.  \nn
\end{eqnarray}
Moreover,
\begin{eqnarray}
\label{eq7.3}
\lambda_{\max} \left( \int_0^T Y_t Y_t' dt \right) = \left\{ \begin{array}{ccc}
O(T^2 ( \ln \ln T)) & \hspace{1cm} {\rm  a.s. \ \ \ if \ \  } \rho =1 \\
O(T^{4\rho -4} (\ln \ln T) ) & \hspace{1cm} {\rm  a.s. \ \ \ if \ \  } \rho \geq 2.
\end{array} \right.
\end{eqnarray}
\end{theorem}

To prove Theorem \ref{th7.1}, we need the following Lemmas.

\begin{lemma}
\label{lm7.1}
\begin{eqnarray}
\sum\limits_{n=j}^\infty \frac{(-1)^n(vt)^{2n-j}}{(2n-j)!} = \left\{ 
\begin{array}{cccc}
O(1) & {\rm  if \ \  } j=0,1 \\ O(t^{j-2}) & {\rm  if \ \  } j\geq 2.
\end{array} \right. \nn
\end{eqnarray}
\end{lemma}

{\bf Proof.}
\begin{eqnarray}
&& \sum\limits_{n=j}^\infty \frac{(-1)^n(vt)^{2n-j}}{(2n-j)!} \nonumber \\
&=& (-1)^j \left[ \frac{(vt)^j}{j!} - \frac{(vt)^{j+2}}{(j+2)!} + \frac{(vt)^{j+4}}{(j+4)!} - \cdots \right] \nonumber \\
&=& \left\{ \begin{array}{cccc}
\cos (vt) & {\rm  if \ \  }j=0 \\
-\sin (vt) & {\rm  if \ \  } j=1 \\
(-1)^{j/2}\left\{ \cos (vt) - \left[ 1- \frac{(vt)^2}{2!} + \cdots +(-1)^{j/2-1}\frac{(vt)^{j-2}}{(j-2)!} \right] \right\}& {\rm  if \ \ } j {\rm \ \ is \
 even, } \ j \geq 2 \\
(-1)^{(j-3)/2}\left\{ \sin (vt) - \left[ vt - \frac{(vt)^3}{3!} + \cdots + (-1)^{(j-1)/2}\frac{(vt)^{j-2}}{(j-2)!} \right] \right\}& {\rm  if \ \ } j {\rm \ \  is \ odd, } \ j \geq 3
\end{array} \right. \nonumber \\
&=& \left\{ \begin{array}{cccc}
 O(1) & {\rm  if \ \  }j=0,1 \\ O(t^{j-2}) & {\rm  if \ \  } j \geq 2.
\end{array} \right. \nn
\end{eqnarray}
Hence, the lemma follows. \qed

\begin{lemma}
\label{lm7.2}
 With the same assumptions as in Theorem \ref{th7.1},
\begin{eqnarray}
||e^{Ft}|| = \left\{ \begin{array}{cccc}
O(1) & {\rm  a.s. \quad if \ \  } \rho =1 \\
O(t^{2\rho -3}) & {\rm  a.s. \quad if \ \  } \rho \geq 2.
\end{array}. \right. \nn
\end{eqnarray}
\end{lemma}

{\bf Proof.}
 Suppose $F$ is a $2\rho \times 2\rho$ matrix and has $\rho$ eigenvalues of
 $\lambda_1=iv$ and $\bar \lambda_1 =-iv$. 
Since the 
characteristic equation for $F$ is $0 = |\la I - F | = (\la - iv)^\rho
 (\la + iv)^\rho = (\la^2 + v^2)^\rho$, by the Cayley-Hamilton theorem,
\begin{eqnarray}
\label{eq7.7}
(F^2 + v^2I)^\rho = 0.
\end{eqnarray}
 {\bf Case 1}: When $\rho =1$, then $F^{2n} = (-1)^n v^{2n}I$ and
\begin{eqnarray}
\label{eq7.8}
e^{Ft} &=& \sum\limits_{n=0}^\infty \frac{F^{2n}t^{2n}}{(2n)!} +
 F\sum\limits_{n=0}^\infty \frac{F^{2n}t^{2n+1}}{(2n+1)!} \nonumber \\
&=& I \sum\limits_{n=0}^\infty 
 \frac{(-1)^n (vt)^{2n}}{(2n)!}
+ \frac{F}{v} \sum\limits_{n=0}^\infty 
 \frac{(-1)^n (vt)^{2n+1}}{(2n+1)!} \nonumber \\
&=& I \cos (vt) + \frac{F}{v} \sin (vt).
\end{eqnarray}
Therefore, $||e^{Ft}|| = O(1)$ when $\rho =1$.  

{\bf Case 2}: When $\rho \geq 2$, then $A = F^2 + v^2I$ is a
 nilpotent matrix of order $\rho$ by (\ref{eq7.7}). Thus,
\begin{eqnarray}
F^2 &=& -v^2 \left[ I -\frac{A}{v^2} \right] \qquad {\rm  and }  \nn\\
F^{2n} &=& (-1)^n v^{2n} \sum\limits_{k=0}^{\rho-1} (-1)^k
 {n\choose{k}}\frac{A^k}{v^{2k}} \nonumber \\
&=& (-1)^n v^{2n}\left( I - \frac{nA}{v^2} + \cdots + (-1)^{\rho-1}
 {n\choose{\rho-1}} \frac{A^{\rho-1}}{v^{2(\rho-1)}} \right). \nn
\end{eqnarray}
Therefore, 
\begin{eqnarray}
\label{eq7.11}
e^{Ft} = \sum\limits_{n=0}^\infty \frac{F^{2n}t^{2n}}{(2n)!} +
 F\sum\limits_{n=0}^\infty \frac{F^{2n}t^{2n+1}}{(2n+1)!}.
\end{eqnarray}

 Let $f_j(n) = 2n (2n-1) \cdots (2n-j+1)$ if $j \geq 1$ and $f_0(n)=1$.
 Then, since \\ $f_0(n), f_1(n), \cdots, f_k(n)$
 are independent, there exist unique $C_0, C_1 \cdots C_k \in {\cal Z}$
 such that
\begin{eqnarray}
{n\choose{k}} = \sum\limits_{j=0}^k C_j f_j(n). \nn
\end{eqnarray}
 Similarly, let $f_j^*(n) = (2n+1) (2n) \cdots (2n-j+2) \ \ {\rm  if } \ \
 j \geq 1$
 and $f_0^*(n)=1$. Then, there exist unique $C_0^*, C_1^*, \cdots C_k^* \in
 {\cal Z}$ such that
\begin{eqnarray}
{n\choose{k}} =\sum\limits_{j=0}^k C_j^* f_j^*(n). \nn
\end{eqnarray}

By Lemma \ref{lm7.1}, the first term of (\ref{eq7.11}) can be expressed as 
\begin{eqnarray}
&&\sum\limits_{n=0}^\infty \frac{(-1)^n (vt)^{2n}}{(2n!)}
 \left[ \sum\limits_{k=0}^{(\rho-1)\wedge n} (-1)^k {n\choose{k}}
 \frac{A^k}{v^{2k}} \right] \nonumber \\
&=& \sum\limits_{k=0}^{\rho-1} \left( -\frac{A}{v^2} \right)^k \left[ \sum\limits_{
n=k}^\infty \frac{(-1)^n (vt)^{2n}}{(2n)!} \left( \sum\limits_{j=0}^k C_j f_j(n)
\right) \right] \hspace{5cm} \nonumber \\ 
&=& \sum\limits_{k=0}^{\rho-1} \left( -\frac{A}{v^2} \right)^k \left[ \sum\limits_{
j=0}^k (vt)^j C_j \left( \sum\limits_{n=k}^\infty \frac{(-1)^n (vt)^{2n-j}}{(2n-
j)!} \right) \right] \nonumber \\
&=& \left\{ \begin{array}{cccc}
 \sum\limits_{k=0}^1 \left( -\frac{A}{v^2} \right)^k \times O(t) &
 {\rm  for \ \ }
 \rho = 2 \\
\sum\limits_{k=0}^{\rho -1} \left( - \frac{A}{v^2} \right)^k \times O(t^{2k-2}) &  {\rm  for \ \ } \rho \geq 3
\end{array} \right. \nonumber \\
&=& \left\{ \begin{array}{cccc}
O(t) & {\rm for \ \ } \rho =2 \\ O(t^{2\rho-4}) & {\rm for \ \ } \rho \geq 3 .
\end{array} \right. \nn
\end{eqnarray}

 Similarly, the second term of (\ref{eq7.11}) can be expressed as
\begin{eqnarray}
&&F\sum\limits_{n=0}^\infty \frac{F^{2n}t^{2n+1}}{(2n+1)!} \nonumber \\
&=& \frac{F}{v} \sum\limits_{n=0}^\infty \frac{(-1)^n (vt)^{2n+1}}{(2n+1)!}
\sum\limits_{k=0}^{(\rho-1)\wedge n } (-1)^k {n\choose{k}} \frac{A^k}{v^{2k}} \hspace{6cm}
 \nonumber \\
&=&
\frac{F}{v} \sum\limits_{k=0}^{\rho-1} \left( -\frac{A}{v^2} \right)^k \left[ 
\sum\limits_{j=0}^k (vt)^j C_j \left( \sum\limits_{n=k}^\infty \frac{(-1)^n
 (vt)^{2n-j+1}}{(2n-j+1)!} \right) \right] \nonumber \\
&=& \frac{F}{v} \sum\limits_{k=0}^{\rho-1} \left( -\frac{A}{v^2} \right)^k \times O(t^{2k-1}) \nonumber \\
&=& O(t^{2\rho-3}). \nn
\end{eqnarray}
Hence, the Lemma follows. 
\qed

\begin{lemma}
\label{lm7.3}
\begin{eqnarray}
\int_0^T (T-s)^k A dW_s =O\left( T^{k+1/2} \sqrt{\ln \ln T}\right) \nn
\end{eqnarray}
\end{lemma}

{\bf Proof.}
 Let $M_u = \int_0^u (t-s)^k AdW_s$, which is a square integrable martingale
 for $[0 < u \leq t]$ and
 $<M>_u = \int_0^u (t-s)^{2k}AA' ds = [t^{2k+1} - (t - u)^{2k+1}] AA'/(2k+1)$.
 Since
 $M_u  =B_{<M>_u}$ by Karatzas and Shreve (\cite{key3} p174),
$$\int_0^T (T-s)^k AdW_s 
= O( B_{T^{2k+1}})
= O(T^{k+1/2} \sqrt{\ln \ln T}).
$$
Hence, the lemma follows. 
\qed

 {\bf Proof of Theorem \ref{th7.1}.}
 If $\rho =1$, then there exist $C \in \cal R$ such that $||e^{Ft}|| \leq C$
 by (\ref{eq7.8}). Therefore,
$$||Y_T||
= ||e^{FT}Y_0 + \int_0^T e^{F(T-s)} AdW_s ||
\leq CY_0 + C\left[ O(\sqrt{T \ln \ln T}) \right] 
= O(\sqrt{T \ln \ln T}).
$$ 
For $\rho \geq 2$, by Lemma \ref{lm7.2} and \ref{lm7.3},
\begin{eqnarray}
||Y_T|| &=& ||e^{FT}Y_0 + \int_0^T e^{F(T-s)}AdW_s || \nonumber \\
& \leq & O\left( ||e^{FT}Y_0||  + ||
\int_0^T \sum\limits_{k=0}^{2\rho -3} C_k (T-s)^k AdW_s || \right) \nonumber \\
&=& O \left( ||e^{FT}Y_0|| + ||
\sum\limits_{k=0}^{2\rho -3} C_k \int_0^T (T-s)^k AdW_s || \right) \nonumber \\
&\leq & O \left( ||e^{FT}Y_0|| + 
\sum\limits_{k=0}^{2\rho -3} |C_k| \times ||O(T^{k+1/2}\sqrt{\ln \ln T})||
 \right) \nonumber \\
&=& O(T^{2\rho -5/2} \sqrt{\ln \ln T}). \nn
\end{eqnarray}

To show (\ref{eq7.3}), we have
\begin{eqnarray}
\lambda_{\max} \left( \int_0^T Y_t Y_t' dt \right) &=& O\left( {\rm tr}
 \int_0^T Y_t Y_t' dt \right) \nonumber \\
&=& O \left( \int_0^T ||Y_t||^2 dt \right) \nonumber \\
&=& \left\{ \begin{array}{cccc}
O(T^2 (\ln \ln T)) & {\rm  a.s. \quad if \ \  } \rho =1 \\
O(T^{4\rho -4} ( \ln \ln T)) & {\rm  a.s. \quad if \ \  } \rho \geq 2.
\end{array} \right. \nn
\end{eqnarray}
Hence, the proof of the theorem. \qed

\subsection{ Zero Eigenvalues}
\vspace{-0.4cm}
In this Section, we study the asymptotic behavior of
the
OU
 processes
when the drift matrix $F$ contains only zeros eigenvalues.(i.e., F is a
 nilpotent matrix.) The main results are summarized in the following:
        
\begin{theorem}
\label{th7.2}
Suppose for the
OU
 process defined in (\ref{eq1.1}),
the RANK condition (\ref{eq1.2}) holds and, all eigenvalues of $F$ are zeros.
 Let $\gamma$ be the
 dimension of the largest block of the rational canonical form of $F$ as
 defined in Section \ref{sec.2} (i.e., $F^\gamma=0$; see the Example). Then
\begin{eqnarray}
||Y_T|| = O(T^{\gamma -1/2} \sqrt{ \ln \ln T}) \qquad {\rm  a.s.} \nn
\end{eqnarray}
Moreover,
\begin{eqnarray}
\label{eq7.22}
\lambda_{\max} \left( \int_0^T Y_t Y_t' dt \right) = O(T^{2\gamma}
 ( \ln \ln T)) \qquad {\rm  a.s.}
\end{eqnarray}
\end{theorem}

{\bf Proof.}
 Since $F$ is a $k \times k$ nilpotent matrix of order $\gamma$
 ($1 \leq \gamma \leq k $), then $F^\gamma =0$ and
\begin{eqnarray}
e^{Ft} =\sum\limits_{n=0}^{\gamma -1} \frac{F^n t^n}{n!} = O(t^{\gamma -1}).
 \nn
\end{eqnarray}
\begin{eqnarray}
||Y_T|| & \leq & O\left( ||e^{FT}Y_0 || + \int_0^T \sum\limits_{k=0}^{\gamma -1}C_k (T-s)^k AdW_s \right) \nonumber \\
&=& O(||e^{FT}Y_0||) + O \left( \sum\limits_{k=0}^{\gamma -1} C_k \int_0^T (T-s)^k AdW_s \right) \nonumber \\
&=& O(T^{\gamma -1}) + O(T^{\gamma -1/2}\sqrt{\ln \ln T}) \nonumber \\
&=& O(T^{\gamma -1/2} \sqrt{\ln \ln T}). \nn
\end{eqnarray}

To prove 
(\ref{eq7.22})
$$\lambda_{\max} \left( \int_0^T Y_t Y_t' dt \right)
= O\left({\rm  Tr} \int_0^T Y_t Y_t' dt \right)
= O\left( \int_0^T ||Y_t||^2 dt \right)
= O(T^{2\gamma} ( \ln \ln T)). 
$$
Hence, the proof. \qed

\end{document}